\documentclass[11pt,reqno]{amsart}
\usepackage{fullpage}
\usepackage{xcolor}
\usepackage[T1]{fontenc}                                   
\usepackage{amsfonts}
\usepackage[utf8]{inputenc}                               
\usepackage{comment}                                       
\usepackage{mparhack}                                      
\usepackage{amsmath,amssymb,amsthm,mathrsfs,eucal}                      

\usepackage{mathtools}
\usepackage{booktabs}                                                                          
\usepackage{graphicx,subfig}                                                                
\usepackage{wrapfig}                                                                            
\usepackage[bookmarks=true,colorlinks=true]{hyperref}                      
\usepackage{bm}

\usepackage{dsfont}
\usepackage{esint}

\newtheorem{definition}{Definition}[section]
\newtheorem{thm}{Theorem}[section]
\newtheorem{prop}{Proposition}[section]
\newtheorem{lem}{Lemma}[section]

\DeclarePairedDelimiter{\abs}{\lvert}{\rvert}
\DeclarePairedDelimiter{\abss}{\bigg\lvert}{\bigg\rvert}
\DeclarePairedDelimiter{\norm}{\lVert}{\rVert}
\newcommand{\R}{\mathbb{R}}
\newcommand{\N}{\mathbb{N}}
\newcommand{\de}{\partial}
\newcommand{\calB}{\mathcal{B}}
\newcommand{\calP}{\mathcal{P}}

\newcommand{\calS}{{\mathsf S}}

\newcommand{\calO}{\mathcal{O}}

\newcommand{\calE}{\mathcal{E}}

\newcommand{\calW}{\mathcal{W}}
\newcommand{\calF}{\mathcal{F}}
\newcommand{\calG}{\mathcal{G}}
\newcommand{\calH}{\mathcal{H}}
\newcommand{\calD}{\mathsf{D}}

\title{A system of continuity equations with nonlocal interactions of Morse type}

\author{Marco Di Francesco, \and Valeria Iorio}

\begin{document}

\address{Marco Di Francesco - DISIM - Department of Information Engineering, Computer Science and Mathematics, University of L'Aquila, Via Vetoio 1 (Coppito)
67100 L'Aquila (AQ) - Italy}
\email{marco.difrancesco@univaq.it}

\address{Valeria Iorio - DISIM - Department of Information Engineering, Computer Science and Mathematics, University of L'Aquila, Via Vetoio 1 (Coppito)
67100 L'Aquila (AQ) - Italy}
\email{valeria.iorio1@univaq.it}



\begin{abstract}
    We study a system of two continuity equations with nonlocal velocity fields using interaction potentials of both attractive and repulsive Morse type. Such a system is of interest in many contexts in multi-population modelling. We prove existence, uniqueness and stability in the $2$-Wasserstein spaces of probability measures via Jordan-Kinderlehrer-Otto scheme and gradient flow solutions in the spirit of the Ambrosio-Gigli-Savaré theory. We then formulate a deterministic particle scheme for this model and prove that gradient flow solutions are obtained in the many particle limit by discrete densities constructed out of moving particles satisfying a suitable system of ODEs. The ODE system is formulated in a non standard way in order to bypass the Lipschitz singularity of the kernel, with difference quotients of the kernel replacing its derivative.
\end{abstract}

\maketitle

\section{Introduction}

Nonlocal aggregation-diffusion equations of the form
\begin{equation}\label{eq:ADE_intro}
    \partial_t \rho =\mathrm{div}(\rho\nabla (a(\rho)+ W\ast \rho))
\end{equation}
arise as a natural modelling tool in many contexts of science and technology in which local and nonlocal interactions among many \emph{particles} or \emph{agents} is considered. Examples of their application may be found in the modelling of swarms \cite{topaz_bertozzi_lewis}, in cellular biology \cite{keller}, in the modelling of crowd dynamics \cite{colombo}, in the physics of granular media \cite{LT}, in material sciences \cite{holm}, in population biology \cite{mogilner}, and in ecology \cite{cantrell}. In some of the above mentioned applications, the use of aggregation-diffusion equations is directly justified to describe a given phenomenon. In other cases, equations of the form \eqref{eq:ADE_intro} may be obtained as suitable scaling limits of more complex models (for example of fluid dynamic type models). Very often these models are recovered via a micro-to-macro limit procedure which is not always rigorously justified, thus motivating a deeper analysis of such a limit.

The mathematical literature on nonlocal aggregation-diffusion equations is extremely rich, and an exhaustive list of references would make this introductory chapter unnecessarily long. We mention here some of the main references for the existence theory, both in a classical functional analytical approach \cite{bertozzi_carrillo_laurent, bertozzi_laurent} and in the context of Wasserstein gradient flows \cite{AGS,CDFLS,CMV}. The formulation and solution of these models without diffusion (that is, with $a=0$ in \eqref{eq:ADE_intro}) is often very helpful to understand the mathematical properties of the nonlocal aggregation part, especially in those cases in which the aggregation kernel $W$ features some singularities. This is quite often the case especially in biological aggregation modelling and in population dynamics. Those singularities often result in the formation of concentrations in finite times \cite{JL}. This makes the analysis of these models more challenging and motivates their study in the framework of measure solutions. We mention here the results in \cite{carrillo_choi_hauray, CDFLS} as examples of well-posedness theories obtained in presence of aggregation kernels featuring various sorts of singularities. 

Among the various aggregation kernels $W$ used in \eqref{eq:ADE_intro} featuring a (mild) singularity, we mention the \emph{Morse kernel}
\begin{equation}\label{eq:Morse_intro}
    W(x)=\pm\frac{1}{2}e^{-|x|},
\end{equation} which is particularly interesting in some applications in swarms dynamics and population biology, see \cite{mogilner}. The plus sign in $W$ models \emph{repulsive} interaction whereas the minus sign models \emph{attractive} interaction. The attractive case is particularly challenging since it yields finite time blow-up of solutions, see \cite{bertozzi_carrillo_laurent,CDFLS}. The repulsive case is part of the theory in \cite{bertozzi_carrillo_laurent} in case of $L^p$ initial data. For the one-dimensional repulsive case in the case of initial data in the space of probability measures we refer to the recent \cite{dis}.

An argument in favour of the use of nonlocal aggregation-diffusion equations is the fact that (at least formally) they feature a \emph{discrete} formulation in the spirit of multi agent systems. More precisely, solutions to equations or systems of the form \eqref{eq:ADE_intro} may often be approximated by time-depending measures or densities constructed out of moving agents, or moving particles, solving a system of ODEs or SDEs. In the case $a=0$, a natural candidate as discrete counterpart of \eqref{eq:ADE_intro} is the ODE system
\begin{equation}\label{eq:ODE_intro}
    \dot{x}_i(t)=-\frac{1}{N}\sum_{j=1}^N \nabla W(x_i(t)-x_j(t))\,,\qquad i=1,\ldots,N\,,
\end{equation}
which describes the movement of $N$ agents with positions $x_i(t)$, $i=1,\ldots,N$. To see this, we observe that \eqref{eq:ODE_intro} may be formally written as
\[\dot{x}_i(t)=-(\nabla W \ast \mu^N(t))(x_i(t))\,,\qquad \mu^N=\frac{1}{N}\sum_{k=1}^N \delta_{x_k(t)}\]
and the above is a natural Lagrangian formulation of the continuity equation \eqref{eq:ADE_intro} in the case $a=0$, with the continuum density $\rho$ replaced by the empirical measure $\mu^N$.

Extensions to \emph{many species} for these kind of systems have been introduced in various contexts, justified by concrete instances occurring in crowd dynamics \cite{degond}, chemotaxis modelling \cite{espejo}, multi-species populations \cite{CM}, and opinion formation \cite{during}. We refer to \cite{DiFEF} for a general mathematical theory of two-species models of the form
\[
\begin{cases}
    \partial_t \rho =\mathrm{div}(\rho\nabla (a(\rho,\eta)+W_{11}\ast \rho + W_{12}\ast \eta)), & \\
    \partial_t \eta =\mathrm{div}(\eta\nabla (b(\rho,\eta)+W_{21}\ast \rho + W_{22}\ast \eta)), &
    \end{cases}
\]
under suitable degenerate parabolicity conditions, see also the earlier \cite{difrafag} for the case without diffusion. A case which was not included in the above theories is the one in which (some of) the kernels $W_{ij}$ feature a repulsive Lipschitz singularity, i.e., they behave like $-|x|$ near zero. Partly inspired by a previous result in the one species case in \cite{BCDP}, the papers \cite{CDEFS,DES} dealt with this case in one space dimension, both for the existence with measure initial condition and with respect to the many particle approximation. 

\medskip
In this paper we contribute to this line of research by considering the case of two species, with the interaction kernels $W_{ij}$ being of Morse type \eqref{eq:Morse_intro}, with a repulsive drift for agents of the same species (repulsive self-interaction) and an attractive one for opposite species (attractive cross-interaction). Hence, we consider the one-dimensional model
\begin{equation}
\label{eq:macroscopic_model}
    \begin{dcases}
        \de_t \rho = \de_x(\rho (W' \ast \rho - W' \ast \eta)), \\
        \de_t \eta = \de_x(\eta (W'\ast \eta -W' \ast \rho)),
    \end{dcases}
\end{equation}
for $(t,x) \in (0, +\infty)\times \R$, with $W$ given by the Morse interaction potential
\begin{equation}
    \label{eq:Morse_potential}
    W(x)=\frac{1}{2} e^{-\abs{x}}\,,
\end{equation}
equipped with the initial datum
\[
\begin{dcases}
    \rho (0,x)=\rho_0 (x), \\
    \eta (0,x)= \eta_0 (x),
\end{dcases}
\]
for $x\in \R$. The unknown is a pair of densities $(\rho(t,x),\eta(t,x))$, modelling two interacting species.
We require $\rho_0, \eta_0$ to be probability measures on $\R$. Additionally, we impose $\rho_0,\eta_0\in L^p(\R)$ for some $p\in (1,+\infty]$. We recall that $W \geq 0$, $W$ is even, and $\int_\R W(x)\,dx = 1$. Moreover, $W$ satisfies the elliptic law
\begin{equation}
    \label{eq:elliptic_law}
    W'' (x) = W(x)-\delta_0 \qquad \mbox{in $\mathcal{D'}(\R)$},
\end{equation}
where $\delta_0$ is the Dirac delta measure centered at $0$. Following the strategy adopted in \cite{CDEFS}, we prove existence and uniqueness of gradient flow solutions in the $2$-Wasserstein space. As usual in this framework, these solutions are constructed by means of the Jordan-Kinderlehrer-Otto (JKO) scheme \cite{JKO}, which is based on a variational formulation of the problem \eqref{eq:macroscopic_model} which uses the functional
\begin{align*}
 \calE [(\rho,\eta)] = &\frac{1}{2} \iint_{\R^2} W (x-y)\,d\rho(x)\,d\rho (y) + \frac{1}{2} \iint_{\R^2} W (x-y) \,d\eta(x)\,d\eta(y) \\
& - \iint_{\R^2} W(x-y)\,d\rho(x)\,d\eta (y)
\end{align*}
and the $2$-Wasserstein metric structure of the space of probability measures with finite second moment. Such a structure results into a \emph{gradient flow} formulation in the spirit of \cite{AGS}. The novel part with respect to the general result of \cite{difrafag} is that the functional $\calE$ does not result from potentials which are \emph{all} convex up to a quadratic perturbation. Indeed, the repulsive cross-interaction terms in the functional result from the repulsive potential $W$ which is not convex up to a quadratic perturbation in view of the decreasing jump of its gradient at the origin. A similar cases was considered in \cite{CDEFS}, but with the Newtonian potential, the analysis of which is made very specific by the fact that the drift of each agent is determined by computing cumulative masses of each species at the agent's position. The case of the Morse potential brings some additional difficulties which motivate a study of its own. We will perform this task and obtain existence, uniqueness, and stability of gradient flow solutions for \eqref{eq:macroscopic_model}.

Then, we analyse the problem of approximating the solutions to \eqref{eq:macroscopic_model} via deterministic interacting particles in the spirit of \eqref{eq:ODE_intro}. In the one species case and with a smooth interaction potential the main reference to this problem is a classical paper by Dobru\v{s}in \cite{dobrusin}. For mildly singular and $\lambda$-convex potentials we refer to \cite{AGS,CDFLS}. The case of the one-dimensional Newtonian potential (both repulsive and attractive) was treated in \cite{BCDP}. We refer to the first part of the recent \cite{dis} for the one species case with the repulsive Morse potential. For the two species case, the result with smooth and mildly singular potentials is part of the results in \cite{difrafag}. The case of singular potentials was treated in \cite{carrillo_choi_hauray}. The recent \cite{DES} studies this problem for a two species system of the form \eqref{eq:macroscopic_model} with $W(x)=-|x|$. 

Following the approach of the aforementioned papers, assuming that both species are made up by $N+1$ particles (but only $N$ of them carry some mass), whose locations are labeled by $x_i$ for the first species and $y_j$ for the second one, as $i,j=0,\ldots, N$,
the (perhaps) most natural deterministic particle approximation of \eqref{eq:macroscopic_model} would be
\begin{equation}
\label{eq:classical_particle_scheme}
    \begin{dcases}
        \dot{x}_i =- \sum_{k\neq i} m_k W' (x_i-x_k) + \sum_k n_k W' (x_i-y_k), \\
        \dot{y}_j =- \sum_{k\neq j} n_k W' (y_j-y_k)+ \sum_k m_k W' (y_j - x_k).
    \end{dcases}
\end{equation}
In the system above, $m_k>0$ are the masses of the particles of the first species and $n_k>0$ are the masses of the particles of the second species. However, in our particle approximation scheme all masses will be required to be equal to $1/N$. In this framework, both species have the same total mass for simplicity.

In \eqref{eq:classical_particle_scheme} the derivative $W$ is not defined at zero and this brings additional difficulties in case one has to deal with particles colliding in a finite time. Therefore, inspired by the recent \cite{dis} for the one-species case (see also the previous \cite{radici_stra}), we will instead consider the alternative particle approximation scheme
\begin{equation}\label{eq:particle_actual}
    \begin{dcases}
        \dot{x}_i = \frac{1}{N} \sum_{k=0}^{N-1} \frac{W (x_{k+1}-x_i) - W(x_k-x_i)}{x_{k+1}-x_k} -\frac{1}{N} \sum_{k=0}^{N-1} \frac{W(y_{k+1}-x_i)-W(y_k - x_i)}{y_{k+1}-y_k}, \\
        \dot{y}_j = \frac{1}{N} \sum_{k=0}^{N-1} \frac{W (y_{k+1}-y_j) - W(y_k-y_j)}{y_{k+1}-y_k} - \frac{1}{N} \sum_{k=0}^{N-1} \frac{W(x_{k+1}-y_j)-W(x_k-y_j)}{x_{k+1}-x_k},
    \end{dcases}
\end{equation}
in which the derivatives of $W$ are replaced by suitably chosen difference quotients. 
This formulation of the particle scheme has the advantage of simplifying the proof of the consistency of the scheme in the many particle $N\rightarrow +\infty$ limit. We remark that the discontinuity of the potential $W$ requires at least some weak $L^p$ estimates on the discrete density in order to prove consistency. We also observe that, unlike other cases (e.g.\ when $W$ is smooth) the many particle limit is not a consequence of the stability in the $2$-Wasserstein distance of gradient flow solutions (which we will prove anyhow) because in this case it is not guaranteed that atomic initial data will produce atomic solutions for all times. The convergence of \eqref{eq:particle_actual} towards \eqref{eq:macroscopic_model} requires uniform estimates on the particle system which are non trivial extensions of what was done in \cite{DES}. Indeed, in the latter the computation of the drift acting on each particle $x_i$ was obtained by simply counting the number of particles of each species on the left and on the right of $x_i$. However, we will be able to get our result by borrowing some of the ideas in \cite{dis} in the one species case.

\medskip
The paper is structured as follows. 
\begin{itemize}
    \item In Section \ref{sec:GF} we analyse the problem of existence, uniqueness, and stability of gradient flow solutions for \eqref{eq:macroscopic_model}. More in detail, in Subsection \ref{subsec:preliminaries} we introduce the main tools we use in optimal transport and gradient flows. and we state the main result of this Section in Theorem \ref{thm:main}. In Subsection \ref{subsec:proof} we prove our main result. The proof uses some convexity property of the functional $\calE$ (Subsubsection \ref{subsubsec:convexity}), the construction via JKO scheme (Subsubsection \ref{subsubsec:JKO}), and the so-called flow-interchange technique introduced in \cite{matthes_mccann_savare} (Subsubsection \ref{subsec:interchange}). The proof is concluded in Subsubsection \ref{subsubsec:conclusion}.
    \item In Section \ref{sec:particles} we analyse the deterministic particle approximation of \eqref{eq:macroscopic_model} via \eqref{eq:particle_actual}. More precisely, in Subsection \ref{subsec:collisions} we prove that no collisions occur for particles of the same species, whereas there may be collisions between particles of opposite species. Subsection \ref{subsec:estimates} is devoted to the main $L^p$ estimate for the scheme \eqref{eq:particle_actual}. Subsection \ref{subsec:limit} is devoted to the proof of Theorem \ref{thm:particles} in which the many particle limit result is stated.
\end{itemize}

\section{Existence and uniqueness of gradient flow solutions}\label{sec:GF}
 
\subsection{Preliminaries and statement of the well-posedness result}\label{subsec:preliminaries}
All concepts defined in this subsection are taken from the book \cite{AGS}. We denote by $\calP(\R^n)$ the space of probability measures on $\R^n$ and by $\calP_2 (\R^n)$ the subspace of probability measures with finite second moment, i.e.,
\[
    \calP_2 (\R^n)= \bigg\{ \mu \in \calP (\R^n) \, : \, \int_{\R^n} \abs{x}^2 \, d\mu(x) < + \infty \bigg\}.
\]
If $\mu \in \calP(\R^n)$ and $T:\R^n \to \R^m$ is a Borel map, we denote by
\[
    \nu = T _\# \mu\,,\qquad \nu\in \calP(\R^m)\,,
\]
the \emph{push-forward measure of $\mu$ through the map $T$}, which is defined as
\[ 
    \nu (A) = \mu (T^{-1} (A)), 
\]
for all Borel sets $A \subset \R^m$. The map $T$ is called a \emph{transport map} pushing the measure $\mu$ to the measure $\nu$. We further restrict to the one-dimensional case for simplicity. The $2$-Wasserstein distance  is defined on $\calP_2(\R)\times \calP_2(\R)$ as
\begin{equation}
\label{2wass}
    W_2 (\mu,\nu)= \bigg( \inf_{\gamma \in \Gamma (\mu,\nu)} \iint_{\R\times\R} \abs{x-y}^2 \,d\gamma (x,y)  \bigg)^{1/2},
\end{equation}
for all $\mu,\nu\in\calP_2(\R)$, where $\Gamma (\mu,\nu)$ is the class of \emph{transport plans} between $\mu$ and $\nu$, i.e.,
\[
    \Gamma (\mu,\nu) = \{ \gamma\in \calP (\R\times\R) \, : \, \pi^1 _\# \gamma = \mu, \calF \pi^2 _\#\gamma=\nu \},
\]
and $\pi^i:\R \times \R \to \R$, $i=1,2$, is the projection operator on the $i$-th component of the product space $\R\times\R$. Denoting by $\Gamma_o (\mu,\nu)$ the class of optimal plans between $\mu$ and $\nu$, namely the minimizers of \eqref{2wass}, the $2$-Wasserstein distance can be rewritten as
\[
    W_2^2 (\mu,\nu) = \iint_{\R \times\R} \abs{x-y}^2 \,d\gamma (x,y),
\]
for all $\gamma \in \Gamma_o (\mu,\nu)$. The existence of optimal plans is guaranteed by Prokhorov's Theorem, see e.g.\ \cite{AGS}.
The pair $(\calP_2 (\R), W_2)$ is a complete metric space. Finally, we denote by $\calP_2^a (\R)$ the set of probability measures with finite second moment that are absolutely continuous with respect to the Lebesgue measure.

Since we are dealing with two interacting species in the one-dimensional space, in this subection we adapt the above definitions to the product space $\calP_2 (\R) ^2$, which we equip with the $2$-Wasserstein product distance defined by
\[
    \calW_2^2 (\bm\mu, \bm \nu ) = W_2^2 (\mu_1, \nu_1) + W_2^2 (\mu_2, \nu_2),
\]
for all $\bm\mu = (\mu_1,\mu_2), \bm \nu=(\nu_1,\nu_2) \in \calP_2 (\R)^2$. Furthermore, given $\bm\mu = (\mu_1, \mu_2)$, we set $L^2 (\bm \mu)= L^2 (d\mu_1) \times L^2 (d\mu_2)$ and
\[
    \norm{\bm v}^2_{L^2 (\bm \mu)} = \int_\R v_1^2(x)\,d\mu_1(x) + \int_\R v_2^2 (x)\,d\mu_2(x), 
\]
for $\bm v = (v_1, v_2) \in L^2 (\bm \mu)$.

Considering $\bm\mu =(\mu_1, \mu_2),\ \bm\nu =(\nu_1, \nu_2) \in \calP_2(\R)^2$, a \emph{constant speed geodesics} connecting $\bm \mu$ and $\bm \nu$ is a curve $\bm \gamma_t:[0,1]\to \calP_2 (\R) ^2$, with $\bm{\gamma}_t = (\gamma_t^1, \gamma_t^2)$ given by
\begin{equation} \label{eq:const_speed_geod}
    \gamma_t^1 = ((1-t)\pi^1+t \pi^2)_\# \gamma_1 \qquad \mbox{and} \qquad \gamma_t^2 = ((1-t)\pi^1+t \pi^2)_\# \gamma_2,
\end{equation}
where $\gamma_1 \in \Gamma_o (\mu_1, \nu_1)$, and $\gamma_2 \in \Gamma_o (\mu_2, \nu_2)$.

\begin{definition}[$\lambda$-convexity along geodesics] Let $\lambda \in \R$, and let $\phi : \calP_2 (\R)^2 \to (-\infty, \infty]$ be a proper functional. We say that $\phi$ is \emph{$\lambda$-geodesically convex} (or simply \emph{$\lambda$-convex}) on $\calP_2 (\R) ^2$ if for any $\bm{\mu}, \bm{\nu} \in \calP_2 (\R) ^2$ there exists a constant speed geodesic $\bm{\gamma}_t$ as in \eqref{eq:const_speed_geod} such that, for all $t\in[0,1]$,
\[
\phi (\bm\gamma_t ) \leq (1-t)\phi (\bm\mu) + t \phi (\bm\nu) - \frac{\lambda}{2} t (1-t)  \calW_2^2 (\bm\mu,\bm\nu).
\]
\end{definition}

\begin{definition}[$k$-flow] A semigroup $S_\phi :[0,\infty]\times\calP_2 (\R) ^2\to \calP_2 (\R)^2$ is a \emph{$k$-flow} for the functional $\phi : \calP_2 (\R)^2 \to (-\infty, \infty]$ with respect to the $2$-Wasserstein distance $\calW_2$ if, for any $\bm\mu \in \calP_2 (\R)^2$, the curve $t\mapsto S_\phi ^t \bm\mu$ is absolutely continuous on $[0,\infty[$ and satisfies the evolution variational inequality (E.V.I.)
\[
    \frac{1}{2} \frac{d^+}{dt} \calW_2^2 (S_\phi^t \bm \mu, \bm \nu) + \frac{k}{2} \calW_2^2 (S_\phi^t \bm \mu, \bm \nu) \leq \phi (\bm \nu) - \phi (S_\phi^t \bm \mu),
\]
for all $t>0$ and for any measure $\bm \nu \in \calP_2 (\R)^2$ with $\phi (\bm \nu) < \infty$.
\end{definition}

Let $\bm \gamma_t \in AC ([0,\infty); \calP_2 (\R)^2)$ be an absolutely continuous curve on $\calP_2 (\R)^2$. Its metric derivative is defined as
\[
\abs{\bm \gamma_t '} (t) \coloneqq \limsup_{h \to 0} \frac{\calW_2 (\bm\gamma_{t+h}, \bm \gamma_t )}{\abs{h}},
\]
and it exists almost everywhere due to the absolutely continuity of $\bm \gamma_t$, see \cite{AGS}.

\begin{definition}[Fréchet sub-differential in $\calP_2(\R)^2$f]
    Let $\phi : \calP_2 (\R)^2 \to (-\infty, \infty]$ be a proper and lower semi-continuous functional. Let $\bm \mu=(\mu_1,\mu_2) \in \calP_2 (\R)^2$. A vector field $ \bm v=(v_1, v_2) \in L^2 (\bm \mu)$ belongs to the \emph{Fréchet sub-differential of $\phi$ at $\bm\mu$}, denoted by $\de \phi (\bm \mu)$, if
\begin{align*}
\phi (\bm \nu) - \phi ( \bm \mu) \geq \inf_{\gamma_i \in \Gamma_o (\mu_i,\nu_i)} & \int_{\R^2 \times\R^2} [v_1(x_1) (y_1-x_1) + v_2 (x_2) (y_2-x_2)]\,d\gamma_1 (x_1,y_1)\,d\gamma_2 (x_2,y_2) \\ 
& + o (\calW_2^2 (\bm \mu,\bm \nu)),
\end{align*}
for all $\bm \nu \in \calP_2 (\R)^2$. Furthermore, if $\de \phi (\bm \mu) \neq \emptyset$, $\de^0 \phi (\bm \mu)$ denotes the element of minimal $L^2 (\bm \mu)$-norm in $\de\phi(\bm \mu)$.
\end{definition}
We observe that $\partial\phi(\bm{\mu})$ is a closed and convex subset of $L^2(\bm{\mu})$, and consequently $\de^0 \phi (\bm \mu)$ is (uniquely) well defined.

If $\phi:\calP_2^a (\R)^2 \to (-\infty, \infty]$ is $\lambda$-convex, for any $\bm\mu \in D(\phi)$ we define the \emph{metric slope} as, cf.\ \cite{AGS},
\[
    \abs{\de \phi} (\bm\mu) = \limsup_{\bm\mu \to \bm\nu} \frac{(\phi (\bm\mu)- \phi (\bm\nu))^+}{\calW_2(\bm\mu,\bm\nu)},
\]
that is finite if and only if $\de \phi (\bm\mu) \neq \emptyset$. Moreover, the metric slope can be rewritten as
\[
    \abs{\de \phi} (\bm\mu) = \min \{ \norm{\bm v}_{L^2 (\bm\mu)} \, : \, \bm v \in \de\phi (\bm\mu) \}\,.
\]
An absolute continuous curve $\bm\gamma (t) : [0,T] \to \calP_2^a (\R)^2$ is a \emph{curve of maximal slope for $\phi$} if the map $t \mapsto \phi (\bm\gamma (t))$ is absolute continuous and it holds
\[
    \phi (\bm\gamma (s)) - \phi (\bm\gamma (t)) \geq \frac{1}{2} \int_s^t \big[ \abs{\bm\gamma '}^2 (\sigma) + \abs{\de \phi}^2 [\bm\gamma (\sigma)] \big]\,d\sigma,
\]
for all $0 \leq s \leq t \leq T$.

We now consider the functional 
\begin{equation} \label{eq:energy_functional}
\begin{aligned}
\calE [(\rho,\eta)] = & \frac{1}{2} \iint_{\R^2} W (x-y)\,d\rho(x)\,d\rho (y) + \frac{1}{2} \iint_{\R^2} W (x-y) \,d\eta(x)\,d\eta(y) \\
& - \iint_{\R^2} W(x-y)\,d\rho(x)\,d\eta (y).
\end{aligned}
\end{equation}
As mentioned in the introduction, \eqref{eq:macroscopic_model} is the formal gradient flow of $\calE$ in the $2$-Wasserstein product space structure. The concept of gradient flows may be formulated in many equivalent ways (under suitable conditions on the functional), including the above defined concept of curve of maximal slope.
Following \cite{AGS}, we formalise this concept in the Definition below.

\begin{definition}[Gradient flow solution]
\label{def:grad_flow_sol}
Let $\bm\gamma_0 = (\rho_0, \eta_0) \in \calP_2^a(\R)^2$. An absolutely continuous curve $\bm\gamma(t,\cdot)=(\rho(t,\cdot),\eta(t,\cdot)):[0,T]\to \calP_2^a(\R)^2$ is a \emph{gradient flow solution to \eqref{eq:macroscopic_model}} if $\rho(t,\cdot)$ and $\eta(t,\cdot)$ solve 
\[
\begin{dcases}
    \de_t \rho (t,x) + \de_x (\rho(t,x)v(t,x)) =0,\\
    \de_t \eta(t,x) + \de_x (\eta(t,x)w(t,x)) =0,
\end{dcases}
\]
in the distributional sense, with initial datum $\bm\gamma_0$, where the velocity field $\bm b(t,x)= (v(t,x), w(t,x))$ satisfies
\[
b_i (t,\cdot)= - (\de^0 \calE [\bm\gamma (t)])_i
\]
for $i=1,2$, and
\[
\norm{\bm b(t,\cdot)}_{L^2(\bm \gamma (t))} = \abs{\bm \gamma '} (t),
\]
for a.e. $t >0$.
\end{definition}

The main result of this section is the following Theorem.

\begin{thm}\label{thm:main}
    Let $T>0$ be fixed and $m \in (1, \infty]$. Assume $(\rho_0, \eta_0) \in (\calP_2 (\R) \cap L^m (\R))^2$. Then, there exists a unique gradient flow solution $(\rho,\eta)\in AC([0,T]; (\calP_2 (\R) \cap L^m (\R) )^2 )$ to system \eqref{eq:macroscopic_model} in the sense of Definition \ref{def:grad_flow_sol}. Moreover, the pair $(\rho, \eta)$ satisfies the properties
    \begin{gather*}
        \int_\R \abs{x}^2 [\rho (t,x) + \eta(t,x)]\,dx  \leq e^{Ct} \int_\R \abs{x}^2 [\rho_0 (x) + \eta_0 (x)]\,dx, \\
        \norm{\rho(t,\cdot)}_{L^m (\R)} + \norm{\eta(t,\cdot)}_{L^m (\R)}  \leq C e^{Ct} \left(\norm{\rho_0}_{L^m (\R)} + \norm{\eta_0}_{L^m (\R)}\right),
    \end{gather*}
    for $t \in [0,T]$, and for some constant $C>0$ independent of $t$. Finally, if $\bm\gamma^1=(\rho^1, \eta^1)$ and $\bm\gamma^2=(\rho^2, \eta^2)$ are two gradient flow solutions to \eqref{eq:macroscopic_model} with initial conditions $\bm \gamma^1_0$ and $\bm \gamma^2_0$ respectively, the following stability estimate holds 
    \begin{equation}\label{eq:stability}
         \calW_2 (\bm\gamma^1(t,\cdot), \bm \gamma^2(t,\cdot)) \leq e^\frac{t}{2} \calW_2 (\gamma^1_0, \gamma^2_0),
    \end{equation}
    for all $t \in [0,T]$.
\end{thm}

\subsection{Proof of Theorem \ref{thm:main}}\label{subsec:proof}
In order to prove existence and uniqueness of solutions to system \eqref{eq:macroscopic_model} in the sense of Definition \ref{def:grad_flow_sol}, we essentially follow the strategy used in \cite{CDEFS} for the case of Newtonian interactions. More precisely,
\begin{itemize}
    \item we start by providing some properties on the interaction energy functional \eqref{eq:energy_functional}, and in particular we characterise the (unique) element of minimal $L^2$-norm of its sub-differential;
    \item we construct our solution as the limit of the so-called \emph{JKO scheme}, see \cite{JKO};
    \item using the so-called \emph{flow interchange technique}, see \cite{matthes_mccann_savare}, we prove some further properties of the limiting curve, in particular a uniform-in-time control on the second order moment and on the $L^p$ norm for an initial condition in $L^p$;
    \item we deduce that the limit of the JKO scheme is a curve of maximal slope, i.e., a gradient flow solution, due to some results in \cite{AGS};
    \item the $\lambda$-convexity property implies a stability estimate in the $2$-Wasserstein space and consequently the uniqueness of solutions.
\end{itemize}

\subsubsection{$\lambda$-convexity of the functional}\label{subsubsec:convexity}
 
\begin{prop} \label{prop:energy_func}
    The functional $\calE$ in \eqref{eq:energy_functional} is $\lambda$-geodesically convex on $\calP_2^a(\R)\times\calP_2^a(\R)$ for all $\lambda\leq -1/2$. Furthermore, setting
\[
    ( \partial^0 W \ast \mu )(x)=\int_{y \neq x} W' (x-y)\,d\mu (y),
\]
for $\mu \in \calP_2(\R)$, the vector field
\begin{equation} \label{eq:element_minimal_norm}
    \partial^0 \calE [(\rho,\eta)] = \begin{pmatrix} \partial^0 W \ast \rho - \partial^0 W \ast \eta \\
    \partial^0 W \ast \eta - \partial^0 W \ast \rho \end{pmatrix}
\end{equation}
is the unique element of minimal $L^2$-norm in the sub-differential of $\calE$.
\end{prop}
\begin{proof}
We write
\[
    \calE [(\rho,\eta)] =  \calH (\rho)+ \calH (\eta)+ \calB [(\rho,\eta)],
\]
with
\begin{align*}
    \calH (\rho) & \coloneqq \frac{1}{2} \iint_{\R \times\R} W(x-y)\,d\rho(x)\,d\rho(y), \\
    \calH (\eta) & \coloneqq \frac{1}{2} \iint_{\R \times\R} W(x-y)\,d\eta(x)\,d\eta (y), \\
    \calB [(\rho,\eta)] & \coloneqq \iint_{\R\times\R} (-W)(x-y)\,d\rho(x)\,d\eta(y).
\end{align*}
    The Morse potential $W(x)$ in \eqref{eq:Morse_potential} can be split as follows:
    \[
    W(x) = S(x)+ N(x)
    \]
    with
    \[
    N(x)= \frac{1}{2} (1-\abs{x}), \qquad S(x)=\frac{1}{2} (e^{-\abs{x}} +\abs{x}-1).
    \]
    The function $N(x)$ (usually referred to as the  Newtonian potential in $1d$) has a Lipschitz singularity at the origin, whereas the function $S(x)$ satisfies $S''(x)=W(x)$ and therefore $S$ belongs to $W^{2,\infty}$. Hence, the functional $\calH$ is $0$-convex seen as a functional on $\calP_2(\R)$, and hence as a functional on $\calP(\R)^2$ as well.
    Concerning the cross term $\calB$, since $-W$ is $-1$-convex, arguing as in \cite[Proposition 4]{CDEFS} (see also \cite[Proposition 3.1]{difrafag}), we get the $-1/2$-geodesic convexity of $\calB$ with $\lambda=-1/2$ because the singular part is $0$-convex (since $-N$ is convex) and the smooth part has second derivative with minimum equal to $-1/2$. Hence, $\calE$ is $-1/2$-convex.
      In order to prove that the right-hand side of   \eqref{eq:element_minimal_norm} belongs to the sub-differential of $\calE$ at $(\rho,\eta)$, by using the additivity of the sub-differential, it is sufficient to show that
\[
    \partial^0 W \ast \rho \in \partial \calH(\rho), \qquad
    \partial^0 W \ast \eta \in \partial \calH (\eta), \qquad
\begin{pmatrix}
    \partial^0 W \ast \eta \\ \partial^0 W \ast \rho
\end{pmatrix} \in \partial \calB [(\rho,\eta)].
\]
Notice that the first two inclusions are meant in the usual one-species sense. These inclusions follow by splitting the Morse potential as $W(x)=N(x)+S(x)$ as above, and using both \cite[Theorem 5.1]{BCDP} for the singular part and the regularity of $S$ and classical results in \cite{AGS} for the smooth part. 
Finally, in order to prove that $\partial^0 \calE$ is the unique element of minimal $L^2$-norm of the sub-differential of $\calE$, one can proceed as in \cite[Proposition 2.2]{CDFLS}, see also \cite[Theorem 10.4.11]{AGS}.
\end{proof}

\subsubsection{JKO scheme}\label{subsubsec:JKO}
Let $\tau >0$ be fixed. Assume that $\bm\gamma_0 = (\rho_0, \eta_0) \in \calP_2^a (\R) \times\calP_2^a (\R)$ is the initial datum and $\calE [\bm\gamma_0] < \infty$. We define recursively the sequence $\{ \bm\gamma_\tau^n \}_{n\in \N} = \{ (\rho_\tau^n, \eta_\tau ^n ) \}_{n\in\N} $ in $\calP_2(\R)^2$ as follows:
\begin{align} 
   & \bm\gamma_\tau ^0 = \bm\gamma_0, \nonumber \\
   & \bm\gamma_\tau ^{n+1} \in \mbox{argmin} \bigg\{ \frac{1}{2\tau} \calW_2^2 (\bm\gamma_\tau^n, \bm\gamma ) + \calE (\bm\gamma), \, \bm\gamma \in \calP_2^a (\R)^2 \bigg\}. \label{eq:jko}
\end{align}
Let $T>0$ and $N = \lceil \frac{T}{\tau} \rceil$. We define the interpolation of the sequence $\{\bm \gamma_\tau^n \}_{n=0,\ldots,N}$ as the piecewise constant curve
\begin{equation}
    \label{eq:gamma_tau}
    \bm\gamma_\tau (t)=\bm\gamma_\tau ^n,
\end{equation}
for $t \in ((n-1)\tau, n \tau]$, as $1\leq n\leq N-1$. We now want to prove that the family of curves $\{\bm \gamma_\tau \}_{\tau >0}$ is compact in a suitable sense. To this aim, we will use a refined version of Ascoli-Arzelà theorem, see \cite[Proposition 3.3.1]{AGS}.

\begin{prop}
    \label{prop:convergence}
    Let $T>0$. The family $\{\bm \gamma_\tau \}_{\tau >0}$ defined in \eqref{eq:gamma_tau} admits a subsequence $\{\bm\gamma_{\tau_k}\}_{k\in \N}$, with $\tau_k\searrow 0$, converging to an absolutely continuous curve $\bm\gamma :[0,T] \to \calP_2 (\R)^2$ uniformly in $t \in [0,T]$ with values in the narrow topology, as $k \to \infty$.
\end{prop}
\begin{proof}
Considering two consecutive iterations $\bm\gamma_\tau^n$ and $\bm\gamma_\tau^{n+1}$, since $\bm\gamma_\tau^{n+1}$ fulfils \eqref{eq:jko}, it holds that
\[
    \frac{1}{2\tau} \calW_2^2 (\bm\gamma_\tau^n, \bm\gamma_\tau^{n+1}) \leq \calE [\bm\gamma_\tau^n] - \calE[\bm\gamma_\tau^{n+1}],
\]
thus
\begin{equation}
    \label{eq:en_bound_in}
    \calE [\bm\gamma_\tau^n] \leq \calE [\bm\gamma_0],
\end{equation}
for all $n \in \N$. Taking $m< n$ and summing over $k$ from $m$ to $n-1$, we get
\[
    \frac{1}{2\tau} \sum_{k=m}^{n-1} \calW_2^2 (\bm\gamma_\tau^k, \bm\gamma_\tau^{k+1}) \leq \calE [\bm\gamma_\tau^m] - \calE [\bm\gamma_\tau ^n].
\]
Considering the definition of $\calE$ and applying H\"older's inequality and Young's inequality for convolutions, we are able to control $\calE [\bm\gamma_t^n]$ from below by means of a constant $\overline{C}$ depending on $\norm{W}_{L^\infty (\R)}$. Therefore, using \eqref{eq:en_bound_in}, we have
\begin{equation}
    \label{eq:sum_bounded}
    \frac{1}{2\tau} \sum_{k=m}^{n-1} \calW_2^2 (\bm\gamma_\tau^k, \bm\gamma_\tau^{k+1}) \leq \calE[\bm\gamma_0] + \overline{C} \eqqcolon \overline{C} (\bm\gamma_0, \norm{W}_{L^\infty (\R)}).
\end{equation}
Now take $m<n$ and $0\leq s < t$, with $s \in ((m-1)\tau, m \tau]$ and $t \in ((n-1)\tau, n\tau]$. By using the inequality \eqref{eq:sum_bounded}, we obtain
\[
    \calW_2^2 (\bm\gamma_\tau (s), \bm\gamma_\tau (t)) = \calW_2^2 (\bm\gamma_\tau^m, \bm\gamma_\tau^n) \leq \bigg[ \sum_{k=m}^{n-1} \calW_2 (\bm\gamma_\tau^k, \bm\gamma_\tau^{k+1}) \bigg]^2 \leq (n-m) \sum_{k=m}^{n-1} \calW_2^2 (\bm\gamma_\tau^k, \bm\gamma_\tau^{k+1}).
\]
For $s=0$, we deduce
\[
    \calW_2^2 (\bm\gamma_0,\bm \gamma_\tau (t)) \leq C(\bm\gamma_0, T, \norm{W}_{L^\infty (\R)}),
\]
that implies that the second moment of $\bm\gamma_\tau (t)$ is uniformly bounded on $[0,T]$, which implies the image of $[0,T]$ via the family of curves $\bm{\gamma}_\tau$ is compact in the narrow topology. Moreover, since $\abs{n-m}< \frac{\abs{t-s}}{\tau} +1$, from \eqref{eq:sum_bounded} we deduce that $\bm\gamma_\tau$ is $1/2$-H\"older equi-continuous, indeed
\[
\calW_2 (\bm\gamma_\tau (s), \bm\gamma_\tau (t)) \leq \abs{n-m}^{1/2} \bigg( \sum_{k=m}^{n-1} \calW_2^2 (\bm\gamma_\tau^k, \bm\gamma_\tau^{k+1})\bigg)^{1/2} \leq c (\sqrt{ \abs{t-s}} + \sqrt{\tau} ),
\]
for some $c>0$. Thus, 
\[
    \limsup_{\tau \to 0^+} \calW_2 (\bm\gamma_\tau (s), \bm\gamma_\tau (t)) \leq \omega (s,t),
\]
where $\omega (s,t) = c \sqrt{\abs{t-s}}$ is a symmetric function on $[0,T]\times [0,T]$, and $\lim_{(s,t)\to (r,r)} \omega (s,t)=0$ for all $r\in [0,T]$.
By applying \cite[Proposition 3.3.1]{AGS}, the statement is proven.
\end{proof}

\subsubsection{Flow interchange}\label{subsec:interchange}
We now want to prove that the curve $\bm\gamma$ obtained in Proposition \ref{prop:convergence} is a gradient flow solution to \eqref{eq:macroscopic_model} in the sense of Definition \ref{def:grad_flow_sol}. To this aim, it is sufficient to show that $\bm\gamma$ is a curve of maximal slope since, from \cite[Theorem 11.1.3]{AGS}, curves of maximal slope coincide with gradient flow solutions in the sense of Definition \ref{def:grad_flow_sol} for $\lambda$-convex functionals. 
We adopt the \emph{flow interchange} strategy, proposed in \cite{matthes_mccann_savare}, that is, we consider some \emph{auxiliary} gradient flows to estimate the dissipation of some other energy functionals we want to control. In particular, we will use this technique twice: first to find a uniform bound in time on the second order moment of $\bm\gamma$, then to obtain $L^m$-regularity for $\bm\gamma$, as $m>1$.

We first consider the decoupled system
\begin{equation}
    \label{eq:flow_inter_2}
    \begin{dcases}
        \de_t u_1 = \de_x (2xu_1),  \\
        \de_t u_2  = \de_x (2xu_2),
    \end{dcases}
\end{equation}
that can be seen as the gradient flow of
\[
    \calG (u_1, u_2)= \int_\R \abs{x}^2 (u_1 (x)+u_2(x))\,dx,
\]
with respect to the $2$-Wasserstein distance $\calW_2$. Denoting by  $\calS^\calG = (\calS^\calG_1, \calS^\calG_2)$ the semigroup generated by system \eqref{eq:flow_inter_2}, we know that $\calS^\calG$ is a $\lambda$-flow for the functional $\calG$ for all $\lambda\geq 0$. We define the dissipation of $\calE$ along $\calS^\calG$ as
\[
    \calD^\calG \calE(\bm\gamma) = \limsup_{h\downarrow 0} \frac{\calE (\bm\gamma) - \calE (\calS^\calG_h \bm\gamma )}{h},
\]
for all $\bm \gamma=(\rho,\eta)\in\calP_2^a (\R) \times\calP_2^a (\R)$.
Consider the continuity equation 
\begin{equation}
\label{eq:ce}
\de_t \zeta - \de_x (2x \zeta )=0,
\end{equation}
with the initial datum $\zeta (t=0)=\zeta_0$. Following \cite[Chapter 8]{AGS}, we can associate to \eqref{eq:ce} to flow map $\Phi_t :\R \to \R$ that is the solution to the ODE model
\[
\label{eq:ce_ode}
\begin{dcases}
    \dot{x} (t)= -2x, \\
    x(t=0)=x_0,
\end{dcases}
\]
i.e., $\Phi_t (x_0)= x_0 e^{-2t}$, and the solution to \eqref{eq:ce} can be represented as
\[
\zeta(t,x)= (\Phi_t) _\# \zeta_0,
\]
namely $\zeta$ is the push-forward of the initial datum through the flow map $\Phi$.

\begin{prop} \label{prop:second_moment}
    Let $T>0$ be fixed. Let $\bm\gamma_0 = (\rho_0, \eta_0) \in\calP_2^a (\R)^2$ be such that $\calG [\bm\gamma_0] < + \infty$. Then, the piecewise constant interpolation $\bm\gamma_\tau = (\rho_\tau, \eta_\tau )$ satisfies
    \[
        \int_\R \abs{x}^2 [ \rho_\tau (t,x) + \eta_\tau (t,x)]\,dx \leq e^{Ct} \int_\R \abs{x}^2 [\rho_0 (x)+\eta_0 (x)]\,dx,
    \]
    for any $t \in [0,T]$, with $C>0$ a constant independent of $t$ and $\tau$. Moreover, the limit $\bm\gamma$ has bounded second order moment uniformly in $[0,T]$.
\end{prop}
\begin{proof}
Since $\bm\gamma_\tau^{n+1}$ is defined as in \eqref{eq:jko}, it holds that
\[
    \frac{1}{2\tau} \calW_2^2 (\bm\gamma_\tau^{n+1} , \bm\gamma_\tau^n ) + \calE (\bm\gamma_\tau ^{n+1}) \leq \frac{1}{2\tau} \calW_2^2 (\calS^\calG_h \bm\gamma_\tau^{n+1} , \bm\gamma_\tau^n ) + \calE (\calS^\calG_h\bm\gamma_\tau ^{n+1}),
\]
for all $h>0$. Considering the definition of dissipation of $\calE$ along $\calS^\calG$, dividing by $h>0$ and taking the $\limsup$ as $h\downarrow 0$, we obtain
\[
    \tau \calD^\calG \calE (\bm\gamma_\tau^{n+1}) \leq \frac{1}{2} \frac{d^+}{dt} \bigg( \calW_2^2 (\calS^\calG_t \bm\gamma_\tau ^{n+1}, \bm\gamma_\tau^n ) \bigg) \bigg\vert_{t=0} \leq \calG (\bm\gamma_\tau^n) - \calG(\bm\gamma_\tau^{n+1}),
\]
where the last inequality follows by the fact that $\calS^\calG$ is a $0$-flow for $\calG$. Concerning the dissipation of $\calE$ along $\calS^\calG$, we get
\begin{equation}
\label{eq:proof_1_dissipation_derivative}
    \calD^\calG \calE (\bm\gamma_\tau^{n+1}) = \limsup_{h\downarrow 0} \frac{\calE (\bm\gamma_\tau ^{n+1}) - \calE (\calS^\calG_h \bm\gamma_\tau^{n+1} )}{h} = \limsup_{h\downarrow 0} \int_0^1 \bigg( - \frac{d}{dz} \bigg\vert_{z=ht} \calE (\calS^\calG_z \bm\gamma_\tau^{n+1}) \bigg)\,dt.
\end{equation}
Estimating the energy functional $\calE$ along the solution $\calS^\calG_t \bm\gamma_\tau^{n+1}$ to \eqref{eq:flow_inter_2}, we get
\begin{align*}
    \calE (\calS^\calG_t \bm\gamma_\tau^{n+1}) & = \calE [(\calS^\calG_{1,t} \rho_\tau ^{n+1} , \calS^\calG_{2,t} \eta_\tau^{n+1})] = \calE [( (\Phi_t)_\# \rho_\tau^{n+1}, (\Phi_t)_\# \eta_\tau^{n+1})] \\
    & = \frac{1}{2} \iint_{\R^2} W (x-y) \,d ((\Phi_t)_\# \rho_\tau^{n+1})(y) \,d ( (\Phi_t)_\# \rho_\tau^{n+1})(x) \\
    & \quad + \frac{1}{2} \iint_{\R^2} W(x-y) \,d ( (\Phi_t)_\# \eta_\tau^{n+1})(x) \,d ((\Phi_t)_\# \eta_\tau^{n+1})(y) \\
    & \quad - \iint_{\R^2} W(x-y) \,d ((\Phi_t)_\# \eta_\tau^{n+1})(y) \,d ( (\Phi_t)_\# \rho_\tau^{n+1})(x) \\
    & = \frac{1}{2} \iint_{\R^2} W(x-y) \,d ( (\Phi_t)_\# \rho_\tau^{n+1}- (\Phi_t)_\# \eta_\tau^{n+1})(y) \,d ( (\Phi_t)_\# \rho_\tau^{n+1}- (\Phi_t)_\# \eta_\tau^{n+1})(x) \\
    & = \frac{1}{2} \iint_{\R^2} W (\Phi_t (x)-\Phi_t (y) ) \,d (\rho_\tau^{n+1} - \eta_\tau^{n+1})(y) \,d (\rho_\tau^{n+1} - \eta_\tau^{n+1})(x) \\
    & = \frac{1}{2} \iint_{\R^2} W ( (x-y) e^{-2t} ) \,d (\rho_\tau^{n+1} - \eta_\tau^{n+1})(y) \,d (\rho_\tau^{n+1} - \eta_\tau^{n+1})(x) \\
    & = \frac{1}{2} \iint_{\R^2} e^{- \abs {x-y} e^{-2t} } \,d (\rho_\tau^{n+1} - \eta_\tau^{n+1})(y) \,d (\rho_\tau^{n+1} - \eta_\tau^{n+1})(x).
\end{align*}
We now have to compute the derivative with respect to time at $t=0$. To this end, we consider the difference quotient and the Taylor expansion of $e^{-\abs{x-y}e^{-2t}}$ in $t=0$ and observe that the ratio
\[
\frac{\calE(\calS^\calG_t \bm\gamma_\tau^{n+1} )-\calE(\bm\gamma_\tau^{n+1} )}{t}
\]
can be passed to the $t\searrow 0$ limit due to Lebesgue dominated convergence's theorem. Hence, we deduce
\[
    \frac{d}{dt} \calE(\calS^\calG_t \bm\gamma_\tau^{n+1} )\bigg\vert_{t=0} = \iint_{\R^2} \abs{x-y} e^{-\abs{x-y}} \,d (\rho_\tau^{n+1} - \eta_\tau^{n+1})(y) \,d (\rho_\tau^{n+1} - \eta_\tau^{n+1})(x).
\]
Since $e^{-\abs{x-y}}\leq 1$ and, by construction, $\bm\gamma_\tau^{n+1} \in \calP_2^a(\R) \times \calP_2^a(\R)$, by H\"older's inequality and Young's inequality, we end up with
\[
\frac{d}{dt} \calE(\calS^\calG_t \bm\gamma_\tau^{n+1} )\bigg\vert_{t=0} \leq C \bigg[ \int_\R \abs{x}^2 \big( \rho_\tau^{n+1} (x) + \eta_\tau^{n+1} (x) \big) \,dx \bigg] = C \calG (\bm\gamma_\tau^{n+1}),
\]
where $C$ is a positive constant. We deduce from \eqref{eq:proof_1_dissipation_derivative}
\[
    (1-C\tau ) \calG (\bm\gamma_\tau^{n+1} ) \leq \calG (\bm\gamma_\tau^n),
\]
and by iteration it holds that
\[
    \calG (\bm\gamma_\tau^n) \leq \bigg( \frac{1}{1-C\tau} \bigg)^n \calG (\bm\gamma_0),
\]
for all $n\in \N$. Since $\tau = T/n$ we obtain for large $n$ (or small $\tau$)
\[
    \bigg( \frac{1}{1-C\tau} \bigg)^n =\bigg[ \bigg( \frac{1}{1-C\frac{T}{n}} \bigg)^\frac{n}{T}\bigg]^T \sim e^{CT},
\]
and thus
\[
    \int_\R \abs{x}^2 [ \rho_\tau^n (x) + \eta_\tau^n (x) ]\,dx \leq e^{CT} \int_\R \abs{x}^2 [ \rho_0 (x) + \eta_0 (x) ]\,dx.
\]
It follows
\[
    \int_\R \abs{x}^2 [ \rho_\tau (t,x) + \eta_\tau (t,x) ]\,dx \leq e^{CT} \int_\R \abs{x}^2 [ \rho_0 (x) + \eta_0 (x) ]\,dx,
\]
for all $t\in [0,T]$, thus the subsequence $\{ \bm\gamma_{\tau_k} \} _{k\in N}$ obtained in Proposition \ref{prop:convergence} has second moment bounded uniformly in $[0,T]$, and this holds also for its limit $\bm\gamma$, due to the weak lower semi-continuity of the second moment w.r.t.\ narrow convergence.
\end{proof}

We now consider the decoupled system 
\begin{equation}
    \label{eq:flow_inter_1}
    \begin{dcases}
        \de_t u_1 = \de_{xx} u_1^m + \varepsilon \de_{xx} u_1, \\
        \de_t u_2 = \de_{xx} u_2^m + \varepsilon \de_{xx} u_2,
    \end{dcases}
\end{equation}
as the gradient flow of the functional
\begin{equation}
    \label{eq:calF}
    \calF (u_1,u_2) = \frac{1}{m-1} \int_\R [u_1 (x)^m + u_2 (x)^m]\,dx + \varepsilon \int_\R [u_1 (x) \log u_1 (x) + u_2 (x) \log u_2 (x)]\,dx,
\end{equation}
with $\varepsilon >0$ and $m\in (1,\infty)$, with respect the $2$-Wasserstein distance $\calW_2$. We denote by $\calS^\calF = (\calS^\calF_1, \calS^\calF_2)$ the semigroup generated by system \eqref{eq:flow_inter_1}, which is well known to be a $0$-flow for the functional $\calF$, see \cite{AGS}. 
For $\bm\gamma=(\rho,\eta)\in\calP_2^a (\R) \times\calP_2^a (\R)$, we define the dissipation of $\calE$ along $\calS^\calF$ as
\[
\calD^\calF \calE (\bm\gamma) = \limsup_{h\downarrow 0} \frac{\calE (\bm\gamma) - \calE (\calS^\calF_h \bm\gamma )}{h}.
\]
In the following proposition we prove that if the initial datum is regular, namely $\bm\gamma_0 \in (\calP_2^a (\R) \cap L^m (\R))^2$, then the piecewise constant interpolation $\bm\gamma_\tau$ keeps this regularity in time.

\begin{lem} \label{lem:entropy_part} 
For an arbitrary $\rho\in \calP_2(\R)$, there exists a positive constant $c>0$ such that
\[
   \int_\R   \rho (\log \rho + c \abs{x}^2)\,dx \geq 0\,.
\]
\end{lem}
\begin{proof}
Let $c$ be a constant we will choose later. For a general $\rho\in \calP_2(\R)$, a straightforward computation yields
\begin{align*}
   \int_\R  & \rho (\log \rho + c \abs{x}^2)\,dx= \int_\R \rho (\log \rho - \log e ^{-c \abs{x}^2} )\,dx \\
    & = \int_\R \frac{\rho}{e ^{-c \abs{x}^2}} \log \frac{\rho}{e ^{-c \abs{x}^2}} e ^{-c \abs{x}^2}\,dx  = \int_\R e ^{-c \abs{x}^2} h \log h \,dx \\
    & = \int_\R \psi (h) e ^{-c \abs{x}^2}\,dx \geq 0,
\end{align*}
since $\psi (h) \coloneqq h \log h - h + 1 \geq 0$, where $h= \rho / e ^{-c \abs{x}^2}$, and the constant $c$ is chosen such that $\int_\R e ^{-c \abs{x}^2}\,dx= 1= \int_\R \rho\,dx$.
\end{proof}

\begin{prop}
    Let $T>0$, and $m \in (1, \infty]$. Assume $\bm\gamma_0= (\rho_0, \eta_0) \in (\calP_2^a (\R) \cap L^m (\R))^2$, with $\calF[\bm\gamma_0] < \infty$. The piecewise constant interpolation $\bm\gamma_\tau = (\rho_\tau, \eta_\tau)$ satisfies
    \[
        \norm{\rho_\tau (t,\cdot)}_{L^m (\R)} +  \norm{\eta_\tau(t,\cdot)}_{L^m (\R)} \leq C e^{Ct} \big(  \norm{\rho_0(\cdot)}_{L^m (\R)} + \norm{\eta_0(\cdot)}_{L^m (\R)} \big),
    \]
    for $m\in (1, \infty)$ and $t\in [0,T]$, where $C$ is a constant depending on $m$. A similar estimate holds for the case $m=+\infty$, with $C=1$.
    Moreover, the limit $\bm\gamma$ belongs to $L^\infty ([0,T], (L^m (\R))^2)$.
\end{prop}
\begin{proof}
Since $\bm\gamma_\tau^{n+1}$ satisfies \eqref{eq:jko}, then
\[
    \frac{1}{2\tau} \calW_2^2 (\bm\gamma_\tau^{n+1} , \bm\gamma_\tau^n ) + \calE (\bm\gamma_\tau ^{n+1}) \leq \frac{1}{2\tau} \calW_2^2 (\calS^\calF_h \bm\gamma_\tau^{n+1} , \bm\gamma_\tau^n ) + \calE (\calS^\calF_h\bm\gamma_\tau ^{n+1}),
\]
for all $h>0$. By the definition of dissipation of $\calE$ along $\calS^\calF$, we get
\begin{equation}
\label{eq:proof_first_inequality}
    \tau \calD^\calF \calE (\bm\gamma_\tau^{n+1}) \leq \frac{1}{2} \frac{d^+}{dt} \bigg( \calW_2^2 (\calS^\calF_t \bm\gamma_\tau ^{n+1}, \bm\gamma_\tau^n ) \bigg) \bigg\vert_{t=0} \leq \calF (\bm\gamma_\tau^n) - \calF(\bm\gamma_\tau^{n+1}),
\end{equation}
where the last inequality follows by the fact that $\calS^\calF$ is a $0$-flow for $\calF$. As in the proof of Proposition \ref{prop:second_moment}, we obtain
\begin{equation}
\label{eq:dissipation_equal_derivative}
    \calD^\calF \calE (\bm\gamma_\tau^{n+1}) = \limsup_{h\downarrow 0} \frac{\calE (\bm\gamma_\tau ^{n+1}) - \calE (\calS^\calF_h \bm\gamma_\tau^{n+1} )}{h} = \limsup_{h\downarrow 0} \int_0^1 \bigg( - \frac{d}{dz} \bigg\vert_{z=ht} \calE (\calS^\calF_z \bm\gamma_\tau^{n+1}) \bigg)\,dt.
\end{equation}
Since $\calS^\calF$ is the flow generated by the solution to \eqref{eq:flow_inter_1}, we compute
\begin{align*}
    \frac{d}{dt} \calE (\calS^\calF_t \bm\gamma_\tau ^{n+1} ) = & \int_\R \de_{xx} [( \calS^\calF_{1,t} \rho_\tau^{n+1} )^m - ( \calS^\calF_{2,t} \eta_\tau^{n+1} )^m ] W \ast ( \calS^\calF_{1,t} \rho_\tau^{n+1} -\calS^\calF_{2,t} \eta_\tau^{n+1} )\,dx \\
    & + \varepsilon \int_\R \de_{xx} ( \calS^\calF_{1,t} \rho_\tau^{n+1} -\calS^\calF_{2,t} \eta_\tau^{n+1} ) W \ast ( \calS^\calF_{1,t} \rho_\tau^{n+1} -\calS^\calF_{2,t} \eta_\tau^{n+1} )\,dx.
\end{align*}
Integrating by parts and using the elliptic law \eqref{eq:elliptic_law}, we get
\begin{align*}
    \frac{d}{dt} \calE (\calS^\calF_t \bm\gamma_\tau ^{n+1} ) = & \int_\R ((\calS^\calF_{1,t} \bm\gamma_\tau^{n+1})^m - (\calS^\calF_{2,t} \bm\gamma_\tau^{n+1})^m ) W \ast (\calS^\calF_{1,t} \bm\gamma_\tau^{n+1} - \calS^\calF_{2,t} \bm\gamma_\tau^{n+1})\,dx \\
    & - \int_\R ((\calS^\calF_{1,t} \bm\gamma_\tau^{n+1})^m - (\calS^\calF_{2,t} \bm\gamma_\tau^{n+1})^m) (\calS^\calF_{1,t} \bm\gamma_\tau^{n+1} - \calS^\calF_{2,t} \bm\gamma_\tau^{n+1})\,dx \\
    & + \varepsilon \int_\R (\calS^\calF_{1,t} \bm\gamma_\tau^{n+1} - \calS^\calF_{2,t} \bm\gamma_\tau^{n+1}) W \ast (\calS^\calF_{1,t} \bm\gamma_\tau^{n+1} - \calS^\calF_{2,t} \bm\gamma_\tau^{n+1})\,dx \\
    & - \varepsilon \int_\R (\calS^\calF_{1,t} \bm\gamma_\tau^{n+1} - \calS^\calF_{2,t} \bm\gamma_\tau^{n+1})^2\,dx \\
    & - \bigg[ ((\calS^\calF_{1,t} \bm\gamma_\tau^{n+1})^m- (\calS^\calF_{2,t} \bm\gamma_\tau^{n+1})^m) W' \ast (\calS^\calF_{1,t} \bm\gamma_\tau^{n+1} - \calS^\calF_{2,t} \bm\gamma_\tau^{n+1}) \bigg]_{x=-\infty}^{x=+\infty} \\
    & + \bigg[ \de_x ((\calS^\calF_{1,t} \bm\gamma_\tau^{n+1})^m - (\calS^\calF_{2,t} \bm\gamma_\tau^{n+1})^m) W \ast (\calS^\calF_{1,t} \bm\gamma_\tau^{n+1}-\calS^\calF_{2,t} \bm\gamma_\tau^{n+1}) \bigg]_{x=-\infty}^{x=+\infty} \\
    & - \varepsilon \bigg[ (\calS^\calF_{1,t} \bm\gamma_\tau^{n+1}-\calS^\calF_{2,t} \bm\gamma_\tau^{n+1})W' \ast (\calS^\calF_{1,t} \bm\gamma_\tau^{n+1}-\calS^\calF_{2,t} \bm\gamma_\tau^{n+1}) \bigg]_{x=-\infty}^{x=+\infty} \\
    & + \varepsilon \bigg[ \de_x (\calS^\calF_{1,t} \bm\gamma_\tau^{n+1}-\calS^\calF_{2,t} \bm\gamma_\tau^{n+1}) W \ast (\calS^\calF_{1,t} \bm\gamma_\tau^{n+1}-\calS^\calF_{2,t} \bm\gamma_\tau^{n+1}) \bigg]_{x=-\infty}^{x=+\infty}.
\end{align*}
The boundary terms vanish since the solution to \eqref{eq:flow_inter_1} decays rapidly at infinity. Indeed
\[
    \abss{ ((\calS^\calF_{1,t} \bm\gamma_\tau^{n+1})^m- (\calS^\calF_{2,t} \bm\gamma_\tau^{n+1})^m) W' \ast (\calS^\calF_{1,t} \bm\gamma_\tau^{n+1} - \calS^\calF_{2,t} \bm\gamma_\tau^{n+1}) } \leq 2 \norm{W'}_{L^\infty (\R)} ((\calS^\calF_{1,t} \bm\gamma_\tau^{n+1})^m + (\calS^\calF_{2,t} \bm\gamma_\tau^{n+1}) ^m)
\]
vanishes as $\abs{x} \to \infty$. Concerning the second boundary term we get
\[
    \abss{ \de_x ((\calS^\calF_{1,t} \bm\gamma_\tau^{n+1})^m - (\calS^\calF_{2,t} \bm\gamma_\tau^{n+1})^m) W \ast (\calS^\calF_{1,t} \bm\gamma_\tau^{n+1}-\calS^\calF_{2,t} \bm\gamma_\tau^{n+1}) } \leq 2 \abs{ \de_x ((\calS^\calF_{1,t} \bm\gamma_\tau^{n+1})^m + (\calS^\calF_{2,t} \bm\gamma_\tau^{n+1})^m) }
\]
that goes to zero as $\abs{x} \to \infty$. The third term can be estimated as
\[
    \abss{ (\calS^\calF_{1,t} \bm\gamma_\tau^{n+1}-\calS^\calF_{2,t} \bm\gamma_\tau^{n+1})W' \ast (\calS^\calF_{1,t} \bm\gamma_\tau^{n+1}-\calS^\calF_{2,t} \bm\gamma_\tau^{n+1}) } \leq 2 \norm{W'}_{L^\infty(\R)} (\calS^\calF_{1,t}\bm\gamma_\tau^{n+1} + \calS^\calF_{2,t} \bm\gamma_\tau^{n+1})
\]
that converges to $0$ as $\abs{x} \to \infty$. About the last term we have
\[
    \abss{ \de_x (\calS^\calF_{1,t} \bm\gamma_\tau^{n+1}-\calS^\calF_{2,t} \bm\gamma_\tau^{n+1}) W \ast (\calS^\calF_{1,t} \bm\gamma_\tau^{n+1}- \calS^\calF_{2,t} \bm\gamma_\tau^{n+1}) } \leq 2 \abs{ \de_x (\calS^\calF_{1,t} \bm\gamma_\tau^{n+1} + \calS^\calF_{2,t} \bm\gamma_\tau^{n+1}) }
\]
that again vanishes as $\abs{x} \to \infty$. Thus, by using H\"older's inequality and Young's inequality for convolutions, and recalling $\norm{W \ast (\calS^\calF_{1,t} \bm\gamma_\tau^{n+1} - \calS^\calF_{2,t} \bm\gamma_\tau^{n+1})} _{L^\infty(\R)} \leq \norm{W}_{L^\infty(\R)} \norm{\calS^\calF_{1,t} \bm\gamma_\tau^{n+1} + \calS^\calF_{2,t} \bm\gamma_\tau^{n+1}}_{L^1(\R)} \leq 1$, we deduce
\begin{align*}
    \frac{d}{dt} \calE (\calS^\calF_t \bm\gamma_\tau ^{n+1} ) \leq & \int_\R ( (\calS^\calF_{1,t} \rho_\tau^{n+1})^m - (\calS^\calF_{2,t} \eta_\tau^{n+1})^m ) W \ast (\calS^\calF_{1,t} \rho_\tau^{n+1} - \calS^\calF_{2,t} \eta_\tau^{n+1})\,dx \\
    & + \varepsilon \int_\R (\calS^\calF_{1,t} \rho_\tau^{n+1} - \calS^\calF_{2,t} \eta_\tau^{n+1})W \ast (\calS^\calF_{1,t} \rho_\tau^{n+1} - \calS^\calF_{2,t} \eta_\tau^{n+1})\,dx  \\
    \leq\ & \int_\R ((\calS^\calF_{1,t} \rho_\tau^{n+1})^m + (\calS^\calF_{2,t} \eta_\tau^{n+1})^m ) \,dx + 2 \varepsilon.
\end{align*}
Thus, we get
\begin{align*}
    \calD^\calF \calE(\bm\gamma_\tau^{n+1}) & \geq \limsup_{h \downarrow 0} \int_0^1 \bigg( - \int_\R [(\calS_{1,z}^\calF \rho_\tau ^{n+1})^m + \calS_{2,z}^\calF \eta_\tau^{n+1})^m]\mid_{z=ht}\,dx\,dt - 2\varepsilon \\
    & = - \liminf_{h\downarrow 0} \int_0^1 \int_\R [(\calS_{1,z}^\calF \rho_\tau ^{n+1})^m + \calS_{2,z}^\calF \eta_\tau^{n+1})^m]\mid_{z=ht}\,dx\,dt - 2\varepsilon.
\end{align*}
By \eqref{eq:proof_first_inequality}-\eqref{eq:dissipation_equal_derivative}, and by using the definition of $\calF$ in \eqref{eq:calF}, we obtain
\[
    (1-\tau (m-1)) \calF (\bm\gamma_\tau ^{n+1}) \leq \calF (\bm\gamma_\tau^n) + 2 \varepsilon \tau - \tau \varepsilon (m-1) \int_\R [\rho_\tau^{n+1}\log \rho_\tau ^{n+1} + \eta_\tau^{n+1}\log \eta_\tau^{n+1}]\,dx.
\]
By adding and subtracting the second moment of both $\rho_\tau^{n+1}$ and $\eta_\tau^{n+1}$ (which are finite at any time due to Proposition \ref{prop:second_moment}), multiplied by a certain constant, from Lemma \ref{lem:entropy_part} we end up with
\[
    (1-\tau (m-1)) \calF (\bm\gamma_\tau ^{n+1}) \leq \calF (\bm\gamma_\tau^n) + C m \varepsilon \tau,
\]
for some constant $C$ depending on $T$, 
which gives
\[
    \calF (\bm\gamma_\tau ^{n+1}) \leq \frac{\calF (\bm\gamma_\tau^n)}{(1-\tau (m-1))} + \frac{C m \varepsilon \tau}{(1-\tau (m-1))}.
\]
Iterating this estimate, we get
\[
    \calF[\bm\gamma_\tau^n] \leq \bigg( \frac{1}{1-\tau (m-1)} \bigg)^n\calF (\bm\gamma_0) + Cm\tau\varepsilon \sum_{k=0}^n \frac{1}{(1-\tau (m-1))^k}.
\]
By using once again the definition of $\calF$, and Lemma \ref{lem:entropy_part}, we obtain (by possibly renaming the constant $C$)
\[
\frac{1}{m-1} \int_\R [(\rho_\tau^n (x))^m + \eta_\tau^n (x))^m]\,dx \leq \bigg( \frac{1}{1-\tau (m-1)} \bigg)^n\calF (\bm\gamma_0) + Cm\tau\varepsilon \sum_{k=0}^n \frac{1}{(1-\tau (m-1))^k} +Cm\varepsilon.
\]
Now, by sending $\varepsilon \searrow 0$ we have
\[
    \frac{1}{m-1} \int_\R [(\rho_\tau^n (x))^m + \eta_\tau^n (x))^m]\,dx \leq \bigg( \frac{1}{1-\tau (m-1)} \bigg)^n \int_\R [(\rho_0 (x))^m + (\eta_0 (x)) ^m]\,dx.
\]
Recalling that $\tau=T/n$, since
\[
    \bigg( \frac{1}{1-\tau (m-1)} \bigg)^n = \bigg[ \bigg( \frac{1}{1-\frac{T(m-1)}{n}} \bigg)^ \frac{n}{T(m-1)} \bigg]^{T(m-1)} \sim e^{T(m-1)},
\]
we get
\[
    \int_\R [ (\rho_\tau^n (x))^m + (\eta_\tau^n (x))^m ]\,dx \leq (m-1) e^{t(m-1)} \int_\R [ (\rho_0 (x))^m + (\eta_0 (x))^m ],
\]
for all $n \in \N$, that implies
\[
    \int_\R [ (\rho_\tau (t,x))^m + (\eta_\tau (t,x))^m ]\,dx \leq (m-1) e^{t(m-1)} \int_\R [ (\rho_0 (x))^m + (\eta_0 (x))^m ],
\]
for $t \in [0,T]$, and $m \in (1,\infty)$. Moreover,
\[
    \bigg( \int_\R [ (\rho_\tau (t,x))^m + (\eta_\tau (t,x))^m ]\,dx \bigg) ^\frac{1}{m} \leq (m-1)^\frac{1}{m} e^{t\frac{m-1}{m}} \bigg( \int_\R [ (\rho_0 (x))^m + (\eta_0 (x))^m ] \bigg)^\frac{1}{m},
\]
from which we deduce that the case $m=+\infty$ is also satisfied by sending $m\rightarrow+\infty$.
We obtain that the subsequence $\{\bm\gamma_{\tau_k}\}_{k\in\N}$ of Proposition \ref{prop:convergence} is bounded in $L^\infty ([0,T], L^m (\R))^2$, thus it admits a converging subsequence $\bm\gamma' \in L^m ([0,T],\R)^2$ in the weak $L^m_{t,x}$ topology in the case $m$ finite. In the case $m=+\infty$, such subsequence exists in the weak-$\ast$ topology of $L^\infty ([0,T] \times \R)$. The limit $\bm\gamma'$ coincides with $\bm\gamma$ on $[0,T]$ due to Proposition \ref{prop:convergence}, and it features the same estimate.
\end{proof}

\subsubsection{Conclusion of the existence and uniqueness proof}\label{subsubsec:conclusion}
In this subsection we want to prove that $\bm\gamma$ is a curve of maximal slope for $\calE$. Let $\Tilde{\bm\gamma}_\tau : [0,+\infty) \to \calP_2(\R)^2$ be the De Giorgi variational interpolation of the discrete values $\{ \bm\gamma_\tau^n\}$ defined in \eqref{eq:jko}, that satisfies
\[
\Tilde{\bm\gamma}_\tau (t)= \Tilde{\bm\gamma}_\tau ((n-1)\tau + \delta) \in \mbox{argmin} \bigg\{\frac{1}{2 \delta} \calW_2^2 (\bm\gamma_\tau ^{n-1}, \bm\gamma) + \calE (\bm\gamma), \, \bm\gamma\in \calP_2 (\R)^2 \bigg\},
\]
if $t= (n-1)\tau + \delta \in ((n-1)\tau, n\tau]$. By arguing as in \cite{AGS}, the following inequality holds:
\begin{equation}
    \label{eq:slope_interpolation}
    \frac{1}{2} \int_0^T \norm{\bm v_{\tau_k} (t)}^2 _{L^2 (\bm\gamma_{\tau_k} (t))} \,dt + \frac{1}{2} \int_0^T \abs{\de \calE}^2 [\Tilde{\bm\gamma}_{\tau_k} (t)] \,dt + \calE (\bm\gamma_{\tau_k} (T)) \leq \calE (\bm\gamma_0),
\end{equation}
where $(\bm\gamma_{\tau_k}, \bm v_{\tau_k})$ is the solution to the continuity equation $\de_t \bm\gamma_{\tau_k} (t) + \mbox{div} (\bm v_{\tau_k}(t) \bm\gamma_k (t))=0$ in the sense of distributions, with $\bm\gamma_{\tau_k}$ from Proposition \ref{prop:convergence}, and $\bm v_{\tau_k}$ the unique velocity field with minimal $L^2(\bm\gamma_{\tau_k} (t))$-norm, see \cite[Theorem 8.3.1, Theorem 8.4.5]{AGS}.  Up to a subsequence, the sequences $\bm\gamma_\tau$ and $\Tilde{\bm\gamma}_\tau$ converge narrowly to the same limit $\bm\gamma$ provided in Proposition \ref{prop:convergence}. Then, by the lower semi-continuity of the slope, see \cite{CDFLS}, one can pass to the limit in \eqref{eq:slope_interpolation} proving that $\bm\gamma$ is a curve of maximal slope. Finally, by \cite[Theorem 11.1.4]{AGS},  we conclude that this curve of maximal slope $\bm\gamma$ is the unique gradient flow solution to \eqref{eq:macroscopic_model} in the sense of Definition \ref{def:grad_flow_sol}. Indeed, if $\bm\gamma^1$ and $\bm\gamma^2$ are two gradient flow solutions to \eqref{eq:macroscopic_model} with initial data $\bm\gamma^1_0$ and $\bm\gamma^2_0$ respectively, by \cite[Theorem 11.1.4]{AGS} we deduce the stability estimate
\eqref{eq:stability}
for all $t\in [0,T]$, and thus the uniqueness of the solution is guaranteed.

\section{Deterministic particle approximation}\label{sec:particles}

In this section we deal with the problem of approximating the gradient flow solution found in Theorem \ref{thm:main} with a discrete density constructed out of a set of particles moving through a system of ordinary differential equations, i.e., a set of \emph{deterministic} particles. As mentioned in the introduction, we assume that all the masses are equal to $1/N$ and we deal with the particle scheme
\begin{equation}
    \label{eq:particle_scheme}
    \begin{dcases}
        \dot{x}_i = \frac{1}{N} \sum_{k=0}^{N-1} \frac{W (x_{k+1}-x_i) - W(x_k-x_i)}{x_{k+1}-x_k} - \sum_{k=0}^{N-1} \frac{1}{N} \frac{W(y_{k+1}-x_i)-W(y_k - x_i)}{y_{k+1}-y_k}, \\
        \dot{y}_j = \frac{1}{N} \sum_{k=0}^{N-1} \frac{W (y_{k+1}-y_j) - W(y_k-y_j)}{y_{k+1}-y_k} - \sum_{k=0}^{N-1} \frac{1}{N} \frac{W(x_{k+1}-y_j)-W(x_k-y_j)}{x_{k+1}-x_k},
    \end{dcases}
\end{equation}
as $i,j=0,\ldots,N$. System \eqref{eq:particle_scheme} is coupled with the initial conditions
\begin{equation}\label{eq:particles_initial}
    x_i(0)=x_{i,0}\,,\quad y_i(0)=y_{i,0}\,,\qquad i=0,\ldots,N.
\end{equation}
In general, the conditions $x_{i,0}\leq x_{i+1,0}$ and $y_{i,0}\leq y_{i+1,0}$ are required for all $i=0,\ldots,N-1$. In our paper we shall always work with the strict inequalities 
\[x_{i,0}< x_{i+1,0}\,,\quad y_{i,0}< 
y_{i+1,0}\,,\qquad i=0,\ldots,N-1\,.\]
We set, for all $i,j\in \{0,\ldots,N-1\}$,
\[ 
d_i = x_{i+1}-x_i, \qquad  r_j = y_{j+1}-y_j, 
\]
and
\[
D_i = \frac{1}{Nd_i} = \frac{1}{N (x_{i+1}-x_i)}, \qquad R_j = \frac{1}{Nr_j} = \frac{1}{N (y_{j+1}-y_j)}.
\]
We then define the piecewise constant densities
\begin{equation}
\label{eq:discrete_densities}
\rho^N (t,x)= \sum_{k=0}^{N-1} D_k (t) \mathbf{1}_{[x_k(t), x_{k+1}(t))}(x), \qquad \eta^N (t,x)= \sum_{k=0}^{N-1} R_k (t) \mathbf{1}_{[y_k(t), y_{k+1}(t))}(x).
\end{equation}
We stress at this stage that $\rho^N(t,\cdot)$ ($\eta^N(t,\cdot)$ respectively) is well defined at a given time $t\geq 0$ if no pairs of particles of the species $x$ ($y$ respectively) are colliding at that time. If so, the two functions $\rho^N(t,\cdot)$ and $\eta^N(t,\cdot)$ are probability measures at every time $t\geq 0$. Assuming for the time being that no collisions occur, we have the expressions
\begin{align*}
    \dot{x}_i & = \frac{1}{N} \sum_{k=0}^{N-1} \frac{W (x_{k+1}-x_i) - W(x_k-x_i)}{x_{k+1}-x_k} - \frac{1}{N} \sum_{k=0}^{N-1} \frac{W (y_{k+1}-x_i)-W(y_k-x_i)}{y_{k+1}-y_k} \\
    & = \sum_{k=0}^{N-1} \int_{x_k}^{x_{k+1}} \frac{1}{N} \frac{1}{d_k} W' (z-x_i)\,dz - \sum_{k=0}^{N-1} \int_{y_k}^{y_{k+1}} \frac{1}{N} \frac{1}{r_k} W' (z-x_i)\,dz \\
    & = \sum_{k=0}^{N-1} \int_{x_k}^{x_{k+1}} D_k W' (z-x_i)\,dz - \sum_{k=0}^{N-1} \int_{y_k}^{y_{k+1}} R_k W' (z-x_i)\,dz \\
    & = \sum_{k=0}^{N-1} \int_\R D_k \mathbf{1}_{[x_k,x_{k+1})} (z) W' (z-x_i)\,dz - \sum_{k=0}^{N-1} \int_\R R_k \mathbf{1}_{[y_k,y_{k+1})} (z) W' (z-x_i)\,dz \\
    & = \int_\R \rho^N (t,z) W' (z-x_i)\,dz - \int_\R \eta^N (t,z) W' (z-x_i)\,dz \\
    & = - W'\ast \rho^N (x_i) + W' \ast \eta^N (x_i).
\end{align*}
In the same way, we find that
\[
\dot{y}_j = - W' \ast \eta^N (y_j)+W' \ast \rho^N (y_j).
\]
Therefore, our particle scheme \eqref{eq:particle_scheme} can be re-written as
\begin{equation}
    \label{eq:scheme_rewritten}
    \begin{dcases}
        \dot{x}_i = - W'\ast \rho^N (x_i) + W' \ast \eta^N (x_i), \\
        \dot{y}_j = - W' \ast \eta^N (y_j)+W' \ast \rho^N (y_j),
    \end{dcases}
\end{equation}
as $i,j=0,\ldots,N$.

\subsection{Analysis of particles collisions}\label{subsec:collisions}
We now assume that particles do not overlap at the initial time and we will prove that particles of the same species do not collide for all $t \geq 0$. Clearly, due to classical Cauchy-Lipschitz theory for ODEs, there exists a time $T\geq 0$ such that \eqref{eq:particle_scheme} (as well as its reformulation \eqref{eq:scheme_rewritten}) has a classical solution. Such a solution can be extended as long as particles do not collide.

\begin{prop}\label{prop:collisions}
    Assume that all the particles are detached at the initial time. Then, particle of the same species do not collide for all $t\geq 0$.
\end{prop}
\begin{proof}
We proceed by contradiction. Let us denote by $x_i$ and $x_{i+1}$ two of the colliding particles, and $t^*>0$ is the collision time, i.e., 
\[
    x_i (t^*)=x_{i+1}(t^*).
\]
We can assume that the first collision between these two particles occurs at time $t^*$, so that in the time interval $(t^*-\varepsilon, t^*)$, for some $\varepsilon >0$, the particles are detached, namely
\[
    x_i (s) < x_{i+1}(s)
\]
for all $s \in (t^*-\varepsilon, t^*)$.

\noindent\textbf{Case 1: particles $x_i$ and $x_{i+1}$ have no particles of the opposite species between them.} Let $y_j$ and $y_{j+1}$ be two consecutive particles of the opposite species such that
\[
   y_j (s) \leq x_i(s) < x_{i+1} (s) \leq y_{j+1} (s)
\]
as $s \in (t^*-\varepsilon, t^*)$. Using the first equation in \eqref{eq:scheme_rewritten}, the fundamental theorem of calculus, and the elliptic equation \eqref{eq:elliptic_law}, we get
\begin{align*}
    \dot{d}_i &  = -[(W'\ast\rho^N)(x_{i+1}(s))-(W'\ast\rho^N)(x_{i}(s))]+[(W'\ast\eta^N)(x_{i+1}(s))-(W'\ast\eta^N)(x_{i}(s))]\\
    & = \int_{x_i}^{x_{i+1}} W'' \ast (-\rho^N (z) + \eta^N (z) )\,dz \\
    & = \int_{x_i}^{x_{i+1}} (\rho^N (z) - \eta^N (z) ) \,dz + \int_{x_i}^{x_{i+1}} W \ast (-\rho^N (z) + \eta^N (z) )\,dz \\
    & \geq \frac{1}{N} - \int_{x_i}^{x_{i+1}} \eta^N (z) \,dz - d_i \norm{W \ast (\rho^N- \eta^N) }_{L^\infty(\R)}\,.
\end{align*}
We stress that the use of the fundamental theorem of calculus is justified by the assumption that particles do not collide, which implies $\rho^N$ and $\eta^N$ are in $L^\infty$ and that makes the convolutions $W''\ast \rho^N$ and $W''\ast \eta^N$ bounded. By using Young's inequality for convolution, we know that 
\[
    \norm{W \ast (\rho^N - \eta^N )}_{L^\infty(\R)} \leq \norm{W}_{L^\infty(\R)} \norm{\rho^N - \eta^N}_{L^1(\R)} \leq 1,
\]
thus
\[
    \dot{d}_i \geq \frac{1}{N} - \int_{x_i}^{x_{i+1}} \eta^N (z) \,dz - d_i.
\]
Now, for $z \in (x_i, x_{i+1})$ and $s\in (t^*-\varepsilon, t^*)$, we have that
\[
    \eta^N (z,s) = \frac{1}{N(y_{j+1}-y_j)} \leq \frac{1}{N(x_{i+1}-x_i)}= \frac{1}{Nd_i}\,.
\]
Therefore
\[
\dot{d}_i \geq \frac{1}{N} - d_i \left[ \frac{1}{N d_i}+1 \right] = -d_i.
\]
By using Gr\"onwall inequality, this implies that 
\[
    d_i (t^*) \geq d_i (t^*-\varepsilon) e^{-\varepsilon},
\]
i.e., $d_i > 0$ at time $t^*$ and this is a contradiction.

\noindent\textbf{Case 2: there is one particle of the opposite species between $x_i$ and $x_{i+1}$ and one between $y_{j-1}$ and $y_{j+1}$ stays far.} Assume that
\[
    y_{j-1}(s)\leq x_i (s) < y_j (s) < x_{i+1} (s) \leq y_{j+1}(s)
\]
for $s \in (t^* -\varepsilon, t^*)$, and $x_i (t^*)=y_j (t^*) = x_{i+1} (t^*)$. Furthermore, we assume that $y_{j-1}$, $y_j$, and $y_{j+1}$ do not collide before the time $t^*$, that is, $y_{j-1}$ and $y_{j+1}$ "remain far" from $x_i$ and $x_{i+1}$ up to the collision time. Thus, there exists a constant $\lambda>0$ such that
\[
    \min\{y_j-y_{j-1},y_{j+1}-y_j\} \geq \lambda
\]
for $s \in (t^*-\varepsilon, t^*)$. Proceeding as in Case 1, we easily obtain
\[
    \dot{d}_i \geq \frac{1}{N} -d_i \left[ \frac{1}{N \lambda} +1 \right]\,.
\]
Hence, a simple ODE argument implies $d_i$ gets close to $\frac{\lambda}{1+N\lambda}$ for large times, which implies it cannot denegerate on finite time. We now assume that only one of the two particles $y_{j-1}$ or $y_{j+1}$ collide with all particles among them at time $t^*$. Assume that particle is $y_{j-1}$. Hence, $y_{j+1}-x_i\geq \lambda$ for some $\lambda>0$.
On $s\in (t^*-\varepsilon,t^*)$ we then have
\begin{align*}
    \dot{d}_i & \geq \frac{1}{N}-d_i-\int_{x_i}^{x_{i+1}}\eta^N(z,s) ds\\
    & \ =  \frac{1}{N}-d_i -\frac{y_j-x_i}{N(y_j-x_i)} - \frac{x_{i+1}-y_j}{N(y_{j+1}-x_i)}\\
    & \ \geq -d_i-\frac{x_{i+1}-y_j}{N\lambda}\geq -\left(1+\frac{1}{N\lambda}\right)d_i\,,
\end{align*}
and the above implies
\[d_i(t^*)\geq d_i(t^*-\varepsilon)e^{-\left(1+\frac{1}{N\lambda}\right)\varepsilon},\]
which is a contradiction.

\noindent\textbf{Case 3. Conclusion.} By possibly interchanging the roles of the two species, we have ruled out the case in which two particles of the $y$ species collide with no $x$-particle in between. Hence, the only possible case left is the one in which, with the notation of Case 2, all particles $y_{j-1}$, $x_i$, $y_j$, $x_{i+1}$, and $y_{j+1}$ collide at the same time. By considering in turns the neighbor 
particles, a simple induction argument implies that a collision between particles of the same species is possible if and only if the two species are alternated and all particles collide at the same time. Indeed, if this is not the case we can always fall in one of the cases considered in Case 2. Assume therefore
\[
    x_0 (s) < y_0 (s) < x_1 (s) < y_1 (s) < \cdots < x_i (s) < y_i(s) < \cdots < x_{N-1}(s) < y_{N-1}(s)<x_N(s)<y_N(s)
\]
for $s \in (t^* - \varepsilon, t^*)$, and $x_i (t^*)  = y_j(t^*)$ for all $i,j = 0,\ldots, N$, i.e., all the particles are detached before $t^*$ and collide at $t^*$. This implies that
\[
\sum_{i=0}^{N-1} d_i(t^*) + \sum_{j=0}^{N-1} r_j(t^*) =0 = x_N(t^*)-x_0(t^*) + y_N(t^*) - y_0(t^*).
\]
Taking the time derivative, we omit the time dependence to keep the notation to a minimum, we obtain
\begin{align*}
    \dot{x}_N -\dot{x}_0 + \dot{y}_N - \dot{y}_0 = & - W' \ast ( \rho^N - \eta^N )(x_N) + W' \ast (\rho^N - \eta^N) (x_0) \\
    & - W' \ast (\eta^N - \rho^N) (y_N) + W' \ast (\eta^N - \rho^N) (y_0) \\
     = & - \int_{x_0}^{x_N} W'' \ast (\rho^N - \eta^N) (z)\,dz - \int_{y_0}^{y_N} W'' \ast (\eta^N - \rho^N)(z)\,dz \\
     = & \int_{x_0}^{x_N} (\rho^N- \eta^N)(z)\,dz + \int_{y_0}^{y_N} (\eta^N-\rho^N)(z)\,dz \\
     & - \int_{x_0}^{x_N} W \ast (\rho^N-\eta^N)(z)\,dz - \int_{y_0}^{y_N} W \ast (\eta^N-\rho^N)(z)\,dz,
\end{align*}
where we used the elliptic law \eqref{eq:elliptic_law}. Now, by using Young's inequality for convolutions, we get
\[
    \dot{x}_N -\dot{x}_0 + \dot{y}_N - \dot{y}_0 \geq 2 - [ x_N - x_0 + y_N - y_0] -  \int_{x_0}^{x_N} \eta^N (z)\,dz - \int_{y_0}^{y_N} \rho^N (z)\,dz.
\]
By considering the particle configuration, we have that $x_N - x_0\leq y_N - x_0$ and by the monotonicity of the integrals, since $\eta^N \geq 0$,
\[
    \int_{x_0}^{x_N} \eta^N (z)\,dz \leq \int_{x_0}^{y_N} \eta^N (z)\,dz = 1,
\]
and similarly
\[
    \int_{y_0}^{x_N} \rho^N (z)\,dz \leq \int_{x_0}^{x_N} \rho^N (z)\,dz = 1.
\]
Thus, we arrive at
\[
    \dot{x}_N -\dot{x}_0 + \dot{y}_N - \dot{y}_0 \geq - [x_N - x_0 + y_N - y_0],
\]
and by Gr\"onwall's inequality we obtain
\[
x_N (t^*)-x_0(t^*) + y_N(t^*) -y_0(t^*) \geq [x_N (t^*-\varepsilon)-x_0(t^*-\varepsilon) + y_N(t^*-\varepsilon) -y_0(t^*-\varepsilon)] e^{-\varepsilon},
\]
that is a contradiction. Hence, no collisions occur in finite times between particle of the same species.
\end{proof}

The result in Proposition \ref{prop:collisions} is of paramount importance since it allows to reformulate the scheme \eqref{eq:particle_scheme} in the form \eqref{eq:scheme_rewritten} with $\rho^N$ and $\eta^N$ defined in \eqref{eq:discrete_densities}. Moreover, such a result also implies global-in-time existence for the unique solution to the ODE system \eqref{eq:particle_scheme} provided particles are initially detached.

Unlike particles of the same species, particles of opposite species may indeed collide in a finite time. To see this, consider the example with two particles of each species, that is $N=1$, with initial condition
\[x_0=-2,\quad x_1=-1,\quad y_0=1,\quad y_1=2\,.\]
We set $f(t)=y_0(t)-x_1(t)$. By symmetry, it is easy to deduce that $f(t)=2y_0(t)$. Hence,
\begin{align*}
    \dot{f}(t) & =2\dot{y}_0(t)=-2W'\ast \eta^1(y_0(t))+2W'\ast\rho^1(y_0(t))\\
    & = -\frac{2}{y_1-y_0}\int_{y_0}^{y_1}W'(y_0-z) dz+\frac{2}{x_1-x_0}\int_{x_0}^{x_1}W'(y_0-z) dz.
\end{align*}
Without restriction, $f(t)\geq 0$ for all $t$ (otherwise, by continuity $f(t_1)=0$ at some time $t_1$ and the proof would be complete). Still due to the symmetry, $d_1=r_1=y_1-y_0=x_1-x_0$ and we can use the decreasing monotonicity of $W$ on the positive half line to get
\begin{align*}
    & \dot{f}(t) = -\frac{2(W(0)-W(d_1))}{d_1}-\frac{2(W(f)-W(f+d_1))}{d_1}\leq -\frac{2(W(0)-W(d_1))}{d_1}\,.
\end{align*}
Moreover, still due to the monotonicity of $W$ and due to its convexity,
\begin{align*}
    \dot{d}_1  & = \dot{y}_1-\dot{y}_0= \frac{2(W(0)-W(d_1))}{d_1}-\frac{W(f+d_1)-W(f+2d_1)}{d_1}+\frac{W(f)-W(f+d_1)}{d_1}\\
    & \leq  \frac{3(W(0)-W(d_1))}{d_1}\,.
\end{align*}
Therefore,
\[\frac{d}{dt}(3f(t)+2d_1(t))\leq 0\,,\]
which implies
\[3 f(t)+2 d_1(t)\leq 3 f(0)+2 d_1(0)\eqqcolon\mu_0\,.\]
In particular, $d_1(t)\leq d_1(0)$ for all $t\geq 0$. Now, we claim that $f(t)$ will vanish at some finite time $t^*$. 
Hence,
\[0\leq 3 f(t)+2 d_1(t)\leq \mu_0\,,\]
which implies $0\leq d_1(t)\leq \mu_0/2$. Recalling
\[\dot{f}(t)\leq-\frac{1-e^{-d_1(t)}}{d_1(t)}\,,\]
we observe that the function $[0,+\infty)\ni d\mapsto g(d)=\frac{1-e^{-d}}{d}$ is continuous, strictly decreasing, with $g(0)=1$ and $g(+\infty)=0$. Hence, on the closed interval $d\in [0,\mu_0/2]$ there holds $g(d)\geq g(\mu_0/2)>0$. This implies
\[\dot{f}(t)\leq - g(\mu_0/2)\,,\]
which implies $f$ becomes zero in a finite time. 

The above example shows that the self-repulsion and cross-attraction forces imply a \emph{mixing} phenomenon which was also observed in \cite{DES} for the case of the Newtonian interaction potential. In short, particles tend to set in alternate species order. 

\subsection{$L^p$ estimates on the particle scheme}\label{subsec:estimates}

Having in mind the particle scheme \eqref{eq:particle_scheme} as a tool to approximate solutions to the PDE system \eqref{eq:macroscopic_model} for large $N$, in this subsection we prove some uniform estimates of the $L^p$ norms of $\rho^N$ and $\eta^N$ with respect to $N$. For future use, we compute
\begin{equation}
    \begin{aligned}
        \label{eq:d_dot}
        \dot{d}_i & = \dot{x}_{i+1}-\dot{x}_i \\
        & = - W' \ast\rho^N (x_{i+1}) + W'\ast \eta^N (x_{i+1}) +W'\ast \rho^N (x_i) - W'\ast \eta^N (x_i) \\
        & = -[W'\ast \rho^N (x_{i+1})-W'\ast \rho^N (x_i)] + [W' \ast \eta^N (x_{i+1})-W' \ast \eta^N (x_i)] \\
        & = \int_{x_i}^{x_{i+1}} \left[- W'' \ast \rho^N (z)+W''\ast \eta^N (z) \right]\,dz,
    \end{aligned}
\end{equation}
and
\begin{equation}
    \begin{aligned}
        \label{eq:D_dot}
        \dot{D}_i = & - \frac{\dot{d}_i}{N d_i^2} \\
        = & - \frac{1}{N d_i^2} \left[ \int_{x_i}^{x_{i+1}} - W'' \ast \rho^N (z)+W''\ast \eta^N (z)\,dz \right] \\
        = & D_i \fint_{x_i}^{x_{i+1}} \left[ W'' \ast \rho^N (z)- W'' \ast \eta^N (z) \right] \,dz.
    \end{aligned}
\end{equation}
Similarly, we have
\begin{equation}
    \label{eq:r_dot}
    \dot{r}_j = \int_{y_j}^{y_{j+1}} [-W'' \ast \eta^N (z)+ W'' \ast \rho^N (z)]\,dz,
\end{equation}
and
\begin{equation}
    \label{eq:R_dot}
    \dot{R}_j = R_j \fint_{y_j}^{y_{j+1}} [W'' \ast \eta^N (z)- W'' \ast \rho^N (z)]\,dz.
\end{equation}

\begin{prop} \label{prop:discrete_Lp_bound}
Let $T >0$. Assume $p \in (1, \infty]$. Then, the approximated density $(\rho^N, \eta^N)$ fulfils
\begin{equation} \label{eq:Lp_discrete_control}
    \norm{\rho^N (t,\cdot)}_{L^p (\R)} + \norm{\eta^N (t,\cdot)}_{L^p (\R)} \leq e^{t} \big[ \norm{\rho^N (0,\cdot)}_{L^p (\R)} + \norm{\eta^N (0,\cdot)}_{L^p (\R)} \big],
\end{equation}
for all $t \in [0, T]$ and for all $N\in \N$.
\end{prop}
\begin{proof}
Let $p \in (1, \infty)$. We use \eqref{eq:d_dot}, \eqref{eq:r_dot}, \eqref{eq:D_dot}, and \eqref{eq:R_dot} above to get
\begin{align*}
     \frac{d}{dt} \int_\R & \left((\rho^N)^p + (\eta^N)^p \right) \,dx \\
     & = \frac{d}{dt} \left[ \sum_{i=0}^{N-1} \int_{x_i}^{x_{i+1}} (D_i)^p \,dx + \sum_{j=0}^{N-1} \int_{y_j}^{y_{j+1}} (R_j)^p \,dx \right] \\
    & = \frac{d}{dt} \left[ \sum_{i=0}^{N-1} D_i^p (x_{i+1}-x_i) + \sum_{j=0}^{N-1} R_j^p (y_{j+1}-y_j) \right] \\
    & = (p-1) \int_\R  {\rho^N (z)}^p [W''\ast\rho^N (z)-W'' \ast \eta^N (z)] \,dz \\
    & \qquad + (p-1) \int_\R  {\eta^N (z)}^p [W'' \ast \eta^N (z)-W'' \ast \rho^N (z)]\,dz \\
    & = (p-1) \int_\R ({\rho^N (z)}^p-{\eta^N (z)}^p) [W''\ast (\rho^N (z)-\eta^N (z))]\,dz \\
    & = (p-1) \int_\R ({\rho^N (z)}^p-{\eta^N (z)}^p) [W \ast (\rho^N (z)-\eta^N (z)) - (\rho^N (z)-\eta^N (z))] \,dz,
\end{align*}
having used the elliptic law \eqref{eq:elliptic_law}. Since 
\[\int_\R({\rho^N (z)}^p-{\eta^N (z)}^p)(\rho^N (z)-\eta^N (z))\,dz\geq 0\,,\]
we use Young's inequality for convolution to estimate
\[\norm{W\ast(\rho^N-\eta^N)}_{L^\infty(\R)}\leq \norm{W}_{L^\infty (\R)}(\norm{\rho^N}_{L^1 (\R)}+\norm{\eta^N}_{L^1 (\R)})=1,\]
and we end up with
\[
\frac{d}{dt} (\norm{\rho^N}_{L^p(\R)}^p + \norm{\eta^N}_{L^p(\R)}^p) \leq (p-1) (\norm{\rho^N}_{L^p(\R)}^p + \norm{\eta^N}_{L^p(\R)}^p)\,.
\]
Then, by applying Gr\"onwall's inequality we obtain
\[
    \left( \int_\R \left( \abs{\rho^N (x)}^p + \abs{\eta^N (x)}^p \right) \right)^\frac{1}{p} \leq e^{\frac{p-1}{p}t} \left( \int_\R \left[\abs{\rho^N (0,x)}^p + \abs{\eta^N (0,x)}^p \right] \,dx, \right)^\frac{1}{p}.
\]
Thus, the estimate \eqref{eq:Lp_discrete_control} holds in case of $p$ finite. The case $p=+\infty$ may be obtained by a standard limiting procedure by letting $p\rightarrow+\infty$.
\end{proof}

\subsection{Many particle limit}\label{subsec:limit}
This subsection is devoted to the proof of the many particle limit. We aim at proving that, given a compactly supported initial condition $(\rho_0,\eta_0)\in (\calP_2(\R)\cap L^p(\R))^2$ for some $p\in (1,+\infty]$, we can approximate the corresponding solution to \eqref{eq:macroscopic_model} by some $(\rho^N,\eta^N)$ of the form \eqref{eq:discrete_densities} constructed out of a set of moving particles $x_i$, $y_i$, $i=0,\ldots,N$ solving \eqref{eq:particle_scheme} as $N\rightarrow+\infty$.

Let $(\rho_0,\eta_0)\in (\calP_2(\R)\cap L^p(\R))^2$ for some $p\in (1,+\infty]$. Assume further that both $\rho_0$ and $\eta_0$ have compact support. As customary in the context of deterministic particle approximations, we atomise the initial conditions as follows. We set
\begin{align*}
    & x_{0,0}=\inf(\mathrm{supp}(\rho_0))\,,\qquad y_{0,0}=\inf(\mathrm{supp}(\eta_0))
\end{align*}
and inductively, for all $i\in \{0,\ldots,N-1\}$,
\begin{align*}
    & x_{i+1,0}=\inf\left\{x\in \R\,:\,\, \int_{x_{i,0}}^x \rho_0(z)dz \geq \frac{1}{N}\right\}\,,\\
    & y_{i+1,0}=\inf\left\{y\in \R\,:\,\, \int_{y_{i,0}}^y \eta_0(z)dz \geq \frac{1}{N}\right\}\,.
\end{align*}
We then consider $x_i(t)$ and $y_i(t)$ as the solutions to \eqref{eq:particle_scheme}-\eqref{eq:particles_initial} for $i=0,\ldots,N$. The result in Proposition \ref{prop:collisions} ensures the particles trajectories exist for all times and that the discrete densities $\rho^N(t,\cdot)$ and $\eta^N(t,\cdot)$ defined in \eqref{eq:discrete_densities} exist for all times $t\geq 0$. We now state our main result in this Section.
\begin{thm}\label{thm:particles}
    Let $T>0$ be fixed. Assume $p \in (1, \infty]$ and $(\rho_0, \eta_0) \in (\calP_2 (\R) \cap L^p (\R))^2$ with compact support. Then, the pair of discrete densities $(\rho^N ,\eta^N)$ converges weakly in $L^p ([0,T]\times\R) ^2$ for $p \in (1, \infty)$, and weakly-$\ast$ in $L^\infty ([0,T]\times\R)^2$, to the unique solution $(\rho,\eta)$ to \eqref{eq:macroscopic_model} in the sense of Definition \ref{def:grad_flow_sol} with initial condition $(\rho_0, \eta_0)$.
\end{thm}

Our strategy extends the procedure in \cite{dis} to the case of two species. Some of the proofs are straightforward generalisations of corresponding results in \cite{dis} and will therefore be omitted. We start with a technical result, the proof of which is similar to the one in \cite[Proposition 3.2]{dis}. We omit the details.
\begin{prop} \label{prop:tech_prop}
Let $\varphi \in C^1 (\R)$ be a test function. Then
\begin{equation} \label{eq:tech_prop_1}
    \frac{d}{dt} \int_\R\varphi(x)\rho^N(t,x) \,dx = - \int_\R \rho^N (t,x) \varphi' (x) W' \ast (\rho^N (t,x)- \eta^N (t,x))\,dx + C_N,
\end{equation}
and
\[
    \frac{d}{dt}  \int_\R\varphi(x)\eta^N(t,x) \,dx = - \int_\R \eta^N (t,x) \varphi' (x) W' \ast (\eta^N (t,x)- \rho^N (t,x))\,dx + C_N,
\]
with
\[
    \abs{C_N} \leq \frac{2 \norm{\varphi'}_{L^\infty (\R)}}{N}.
\]
\end{prop}
The result in Proposition \ref{prop:tech_prop} is the basic tool to obtain consistency in the limit, that is, to obtain the pair of continuity equations in \eqref{eq:macroscopic_model} in the $N\rightarrow+\infty$ limit.

The approximation result for the initial datum is stated in the next two Lemmas. For the proofs, see \cite[Lemma 4.1, Lemma 4.2]{dis}.

\begin{lem} \label{lem:discrete_bound_i_d}
    Let $p \in [1, \infty]$. Then the $L^p$-norms of $\rho^N(\cdot,0)$ and $\eta^N(\cdot,0)$ are uniformly bounded by a positive constant $C$ depending on $\norm{\rho_0}_{L^p(\R)}$ and $\norm{\eta_0}_{L^p(\R)}$, i.e.,
    \[
        \norm{\rho^N_0}_{L^p (\R)} \leq C, \qquad \norm{\eta^N_0}_{L^p (\R)} \leq C,
    \]
    for all $N \in \N$.
\end{lem}

\begin{lem} \label{lem:discrete_conv_i_d}
    Let $p\in (1,\infty]$ and assume $(\rho_0, \eta_0) \in (\calP_2 (\R) \cap L^p(\R))^2$. Then, $\rho^N_0 \to \rho_0$ and $\eta^N_0\to \eta_0$ in the sense of distributions. Furthermore, $\rho^N_0 \rightharpoonup \rho_0$ and $\eta^N_0 \rightharpoonup \eta_0$ weakly in $L^p(\R)$ if $p \in (1,\infty)$, whereas $\rho^N_0 \xrightharpoonup{\ast} \rho_0$ and $\eta^N_0 \xrightharpoonup{\ast} \eta_0$ weakly-$\ast$ in $L^\infty (\R)$.
\end{lem}

We are now ready to prove Theorem \ref{thm:particles}.

\begin{proof}[Proof of Theorem \ref{thm:particles}]
    Concerning the initial data, from Lemma \ref{lem:discrete_conv_i_d} we know that
    \[
        \rho^N_0 \rightharpoonup \rho_0, \qquad \mbox{and} \qquad \eta^N_0 \rightharpoonup \eta_0,
    \]
    weakly in $L^p(\R)$ up to a subsequence.
    By Proposition \ref{prop:discrete_Lp_bound} and Lemma \ref{lem:discrete_bound_i_d}, we deduce that there exist two subsequences of both $\rho^N$ and $\eta^N$ converging weakly (or weakly-$\ast$ if $p=+\infty$) to some limits $\rho, \eta \in L^p ([0,T]\times\R)$ as $N\to \infty$. Hence, still denoting by $\rho^N$ and $\eta^N$ such subsequences, we get
    \begin{gather*}
        \lim_{N\to\infty} \abs{W' \ast \rho^N - W' \ast\rho} = \lim_{N\to \infty} \abs{W' \ast (\rho^N - \rho)}=0, \\
        \lim_{N\to\infty} \abs{W' \ast \eta^N - W' \ast\eta} = \lim_{N\to \infty} \abs{W' \ast (\eta^N - \eta)}=0, 
    \end{gather*}
    due to $W'\in L^q (\R)$, where $\frac{1}{p}+\frac{1}{q} =1$. Since by Young's inequality for convolution it holds that 
    \[
        \norm{W' \ast\rho^N}_{L^\infty (\R)} \leq \norm{W'}_{L^q(\R)}, \qquad  \mbox{and} \qquad \norm{W' \ast\eta^N}_{L^\infty (\R)} \leq \norm{W'}_{L^q(\R)},
    \]
    we get
    \[
        W' \ast\rho^N \to W' \ast \rho, \qquad \mbox{and} \qquad W' \ast\eta^N \to W' \ast \eta,
    \]
    in $L^q_{\text{loc}}(\R)$. Considering the result in Proposition \ref{prop:tech_prop}, multiplying \eqref{eq:tech_prop_1} by $\chi=\chi (t) \in C_c^\infty ([0,\infty))$ and integrating on $[0,T]$ with $\text{supp}(\chi) \subset [0,T]$, we end up with
    \[
    \int_0^T \chi(t) \frac{d}{dt} \int_\R \varphi(x) \rho^N \,dx\,dt= - \int_0^T \int_\R \chi(t) \rho^N \varphi' W' \ast (\rho^N - \eta^N) \,dx\,dt + \calO (N^{-1}T).
    \]
    Taking $\psi (t,x)=\chi(t)\varphi(x)$, integrating by parts we obtain
    \[
    \int_0^T \int_\R [ \partial_t \psi- \partial_x \psi W' \ast (\rho^N - \eta^N)] \rho^N \,dx\,dt + \int_\R \psi (0,x)\rho^N_0\,dx = \calO (N^{-1}T).
    \]
    Passing to the limit as $N \to \infty$, we get the desired convergence result for a general test function via standard density of cylindrical test functions. By weak-lower semicontinuity, the limit $(\rho,\eta)$ belongs to $(L^p([0,T]\times\R))^2$ and it therefore coincides with the unique gradient flow solution according to Definition \ref{def:grad_flow_sol}.
\end{proof}

\section*{Acknowledgments}
MDF is partially supported by the Italian “National Centre for HPC, Big Data and Quantum Computing” - Spoke 5 “Environment and Natural Disasters” and by the Ministry of University and Research (MIUR) of Italy under the grant PRIN 2020- Project N. 20204NT8W4, Nonlinear Evolutions PDEs, fluid
dynamics and transport equations: theoretical foundations and applications. VI is supported by the ``MMEAN-FIELDSS'' INdAM project N.E53C22001930001. 
This research is also partially supported by the InterMaths Network, \url{www.intermaths.eu}


\end{document}